\documentclass[12pt]{article}
\usepackage[applemac]{inputenc}
\usepackage{lscape}
\usepackage[dvips]{graphicx}
\usepackage[psamsfonts]{amssymb}
\usepackage[all]{xy}

\newcommand{\bk}{l\!k}

\title{String topology on Gorenstein spaces}
\author{Yves F\'elix and Jean-Claude Thomas}

\begin{document}
\maketitle

\begin{abstract}   The purpose of this paper  is to describe a
general and simple setting for defining $(g,p+q)$-string
operations on a Poincar\'e duality space and more generally on a
Gorenstein space. Gorenstein spaces include
 Poincar\'e duality spaces as well as classifying spaces or homotopy quotients
  of  connected Lie groups.
   Our presentation implies directly  the homotopy invariance of each
    $(g,p+q)$-string operation   as well as  it leads to
     explicit computations.\end{abstract}

Shriek maps play a central role in string topology and its
generalizations. Following  Dold, various presentations have been given
(\cite{Br},\cite{C}, \cite{CK}). Usually shriek maps are defined in
(co)homology for maps $f : N \to M$  from a closed oriented
$n$-manifold to a closed oriented $m$-manifold. Here we will consider shriek maps
at the   cochain level and in a more general setting.

More precisely we   work  in the category   of (left or right)
differential graded modules over a differential graded
$\bk$-algebra $(R,d)$ ($\bk$ is a fixed field), that we call for
sake of simplicity the category of $(R,d)$-modules.   Its
associated derived category is obtained by formally inverting
quasi-isomorphisms, i.e. the maps that induce isomorphisms in
homology and hereafter denoted by $\simeq$ in the diagrams. In the
derived category  the vector space  of  homotopy classes of maps
of degree $q$
   from the $(R,d)$-module $(P,d)$ to the $(R,d) $-module $(Q,d)$ is then denoted
by $\mbox{Ext}^q_R( P,Q)$. In other words an element   $  \varphi
\in \mbox{Ext}^q_R(P,Q)$ is represented by a morphism of
differential $R$-modules $P' \to Q$ (also denoted $\varphi$ by
abuse) where $(P',d)$ is some cofibrant replacement  of $(P,d)$.
The induced map $H(P,d)\cong H(P',d)\to H(Q,d)$ will be denoted
$H(\varphi)$.

For recall,   an {\it oriented Poincar\'e duality space of (formal)
dimension $m$}  is a path connected space  $M$ together with an
{\it orientation class} $[M] \in H_m(M)$ such that the cap product
by the orientation class,
$$
-\cap [M] : H^\ast(M) \to H_{\ast-m}(M) \,,
$$
is an isomorphism.  The {\it fundamental class of $M$} is the
element $ \omega_M \in H^m(M) $ such that   $ <\omega_M, [M]>= 1$
where $<\,\,,\,> :H^\ast(M)\otimes H_\ast(M)\to \bk$ denotes the
Kronecker product. The definition of  an oriented homotopy type is
clear from the above definition.

We consider a pullback diagram,

\centerline{ \footnotesize  $
\xymatrix{
 \ar@{}[d]^{(\ast)} &E'   \ar[d]_{p'}  \ar[rr]^g &&E \ar[d]^p&\ar@{}[d]^{,} \\
 &N\ar[rr]^f &&M&}
$}
where
\centerline{$
(H)\hspace{5mm}\left\{
\begin{array}{l}
  N \mbox{ is an oriented  Poincar\'e duality space of dimension } n\,,\\
 M  \mbox{ is a  1-connected   oriented  Poincar\'e duality space of dimension } m\,,\\
p: E\to M \mbox{ is a fibration}\,,\\
\mbox{$H^*(E)$ is a graded  vector space of finite type}.
\end{array}\right.
$}

\vspace{2mm} \noindent{\bf Theorem A.} {\sl With the notation
above there exist unique elements
$$
f^! \in \mbox{\rm Ext}^{n-m} _{C^\ast(M)}(C^\ast(N), C^\ast(M))
\mbox{ and } g^! \in\mbox{\rm  Ext}^{n-m}_{C^\ast(E)}(C^\ast(E'),
C^\ast(E))
$$
satisfying:
\begin{enumerate}
\item $H^*(f^!)(\omega_N)= \omega_M$;
\item  In  the derived category of $C^\ast(M)$-modules the following
diagram commutes,

\centerline{\footnotesize $
\xymatrix{
  C^\ast(E')    \ar[rr]^{g^!} &&C^{*+n-m} (E )&\ar@{}[d]^{.} \\
 C^\ast(N) \ar[u]^{C^\ast(p')} \ar[rr]^{f^!} &&C^{*+n-m}(M){\,.}
\ar[u]_{C^\ast(p)}& }
$}

\end{enumerate}
Moreover, if the homotopy fiber of the map $f$ is   a Poincar\'e
duality space then $H^*(f^!)$ and $H^*(g^!)$ coincide  with the
respective integrations along the fiber.}

\vspace{2mm}  It is immediate from  theorem A that  if $M$ and $N$
are    closed oriented manifolds  then
  $H^*(f^!)$ is the usual cohomology shriek map in the sense
  of \cite[Chap.V$\!$I-Def.11.2]{Br}.

 \vspace{2mm} As a first example, consider   a smooth embedding of compact oriented manifolds
  $f : N^n \to M^m$. Let  $ V$  be a tubular
neighborhood of $f(N)$ in $M$. Denote    by $\Omega \in C^{m-n}(V,
\partial V)$ a cocycle representing the Thom class of the relative bundle
$(V, \partial V) \to N$. Then we have a sequence of maps of
$C^*(M)$-modules, the second one being the Thom quasi-isomorphism,
$$ C^*(N)\stackrel{\simeq}{\leftarrow} C^*(V) \stackrel{\cup \Omega }{\to}
 C^{*+m-n}(V, \partial V)
\stackrel{\simeq}{\leftarrow} C^{*+m-n}(M, M\smallsetminus f(N))
\hookrightarrow C^{*+m-n}(M)\,.
$$
This sequence induces a $H^\ast(M)$-linear map  $H^\ast(N) \to
H^{\ast+(m-n)} (M) $ which maps the fundamental class to the
fundamental class. Therefore this sequence  defines in the derived
category a representative of the shriek map $f^!$. Another
representative of the same class
  can be
constructed as follows. Let $\psi : (P,d) \stackrel{\simeq}{\to}
C^*(N)$ be a cofibrant replacement of $C^*(N)$ as a
$C^*(M)$-module.  We denote by cap$_M$ and cap$_N$ the cap
products with  representatives of the orientation classes in
$C_*(M)$ and $C_*(N)$.  The lifting property of cofibrant models
(see §1.2-(SF3))  furnishes the following homotopy commutative
diagram of $C_*(M)$-modules,

$$
\xymatrix{
(P,d)\ar[d]^{\simeq}_{\psi} \ar@{.>}[rr]^{\varphi}
&& C^{*+m-n}(M)\ar[d]_{\simeq}^{\mbox{\tiny cap} _M}&\ar@{}[d]^{.}\\
C^*(N)\ar[r]_{\simeq}^{\mbox{\tiny cap}_N } & C_{n-*}(N)\ar[r]^{C_\ast(f)}& C_{n-*}(M)\,.&\\
}
$$
 By construction, $\varphi$ is unique up to
homotopy. Moreover, $H^\ast(\varphi): H^*(N) \to H^*(M)$ maps the
fundamental class to the fundamental class. Therefore $\varphi $
is also a representative  of $f^!$.

As a second example,   consider  a fibration  $ p: E \to M$ with
base  a 1-connected oriented Poincar\'e duality space of dimension
$m$ and  fibre $ F:=p^{-1}(\{b_0\})$. Applying  Theorem A to the
pullback diagram
$$
\xymatrix{
 \ar@{}[d]^{(\dagger)} &F   \ar[d]   \ar[rr]^{j} &&E \ar[d]^p \\
 &\{b_0\}\ar@{^(->}[rr] &&M }
$$
  yields a well defined  homomorphism
$H^*(j^!): H^*(F)\to H^{*+m}(E)$   called {\it the intersection
map with the fiber}.  This  morphism generalizes the intersection
map defined by Chas and Sullivan (\cite{CS}). See \cite{K2} for
another description.

\vspace{2mm} Theorem A can be used to give a homotopy invariant
definition of   {\it the loop product}  on a 1-connected oriented
Poincar\'e duality space $M$ of dimension $m$. Consider the
pullback diagram,
$$\xymatrix{
\ar@{}[d]_{(\ast\ast)}&LM\times_MLM  \ar[d]_{p'} \ar[r]^{q}  & LM\times LM \ar[d]^{p\times p}&\ar@{}[d]^{,}\\
&M  \ar[r]^{\Delta}& M\times M& }$$
 where $\Delta : M\to M\times
M$ denotes the diagonal map, $LM$ is the free loop space, $LM=
M^{S^1}$, and $p : LM\to M$ is the usual fibration that associates
to a loop its base point. Diagram $(\ast\ast)$ sastifies
conditions (H) above. Therefore, by Theorem A,  there exists a
unique  class
$$
{q}^! \in \mbox{Ext}^m_{C^*(LM\times LM)} \left(C^*(LM\times_MLM
), C^*(LM\times LM) \right)\,.
$$
 Denote by $c:   LM\times_MLM  \to LM$ the   composition of free
loops. Then the  linear map of degree $m$
$$ H^*({q}^!)\circ H^*(c) :  H^*(LM)  \stackrel{ }{\longrightarrow}
H^*(LM)\otimes H^*(LM) $$  coincides with   the   dual of the
  loop product   (coefficients   in the field $\bk$)
$$
\bullet :  H_*(LM)\otimes H_\ast(LM) \to H_{*-m}(LM)
  $$
\noindent when $M$ is closed oriented manifold (\cite{CS},
\cite{CJY}, \cite{FT2}).

More generally,  let $S$ be a  connected surface of genus $g$ with $p+q$
boundary components considered as a cobordism between the union of
$p\geq 1$ incoming circles   and   $q\geq 1$ outgoing circles. The embedding
of the incoming circles leads to a cofibration
 $j :\amalg_p S^1 \to S$
inducing for each space $M$ a fibration $q_S : \mbox{Map}(S,M) \to (LM)^p\,.$

For instance,  for $g=2$, $p=3$ and $q=2$,

 \centerline{
 \includegraphics[width=70mm]{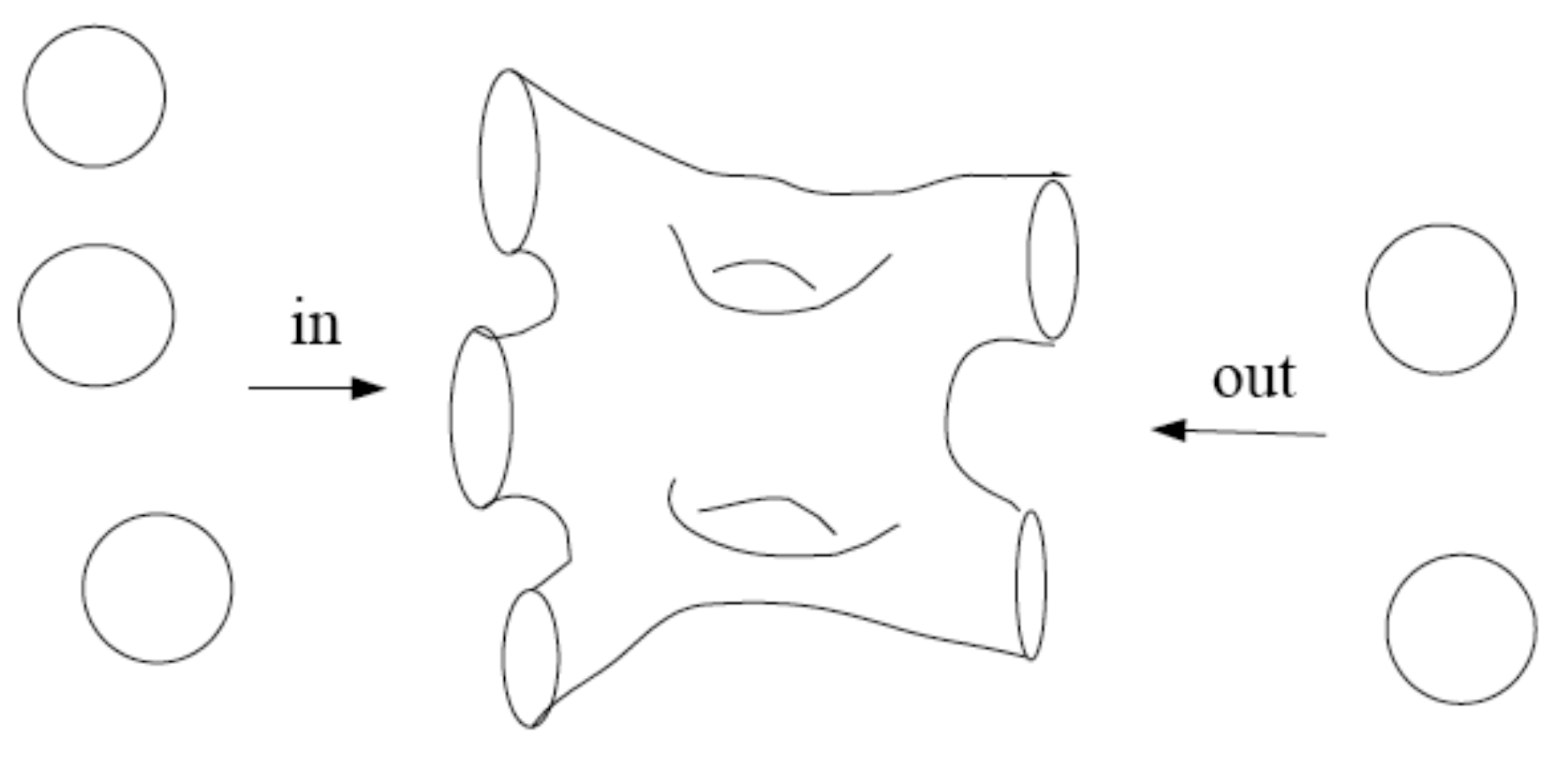}
 }

   Using Sullivan chord diagrams (\cite{CG}) the fibration
$q_S$ can be viewed as a pullback fibration in a diagram of
fibrations
$$\xymatrix{
\mbox{Map}(S,M)\ar[d]_{\scriptstyle q_S}\ar[r] &    M^r
\ar[d]_{\scriptstyle \Delta}\\
 (LM)^p\ar[r]^{\psi}  &M^t }
  $$
 Recall that a Sullivan chord diagram   of the surface  $S$ is
a graph, $C$, homotopy equivalent to $S$ and composed of $p$
circles and $r$ trees $T_i$ whose extremities are $t$ distinct
points on the circles.
   For instance in the chord diagram

 \centerline{
 \includegraphics[width=50mm]{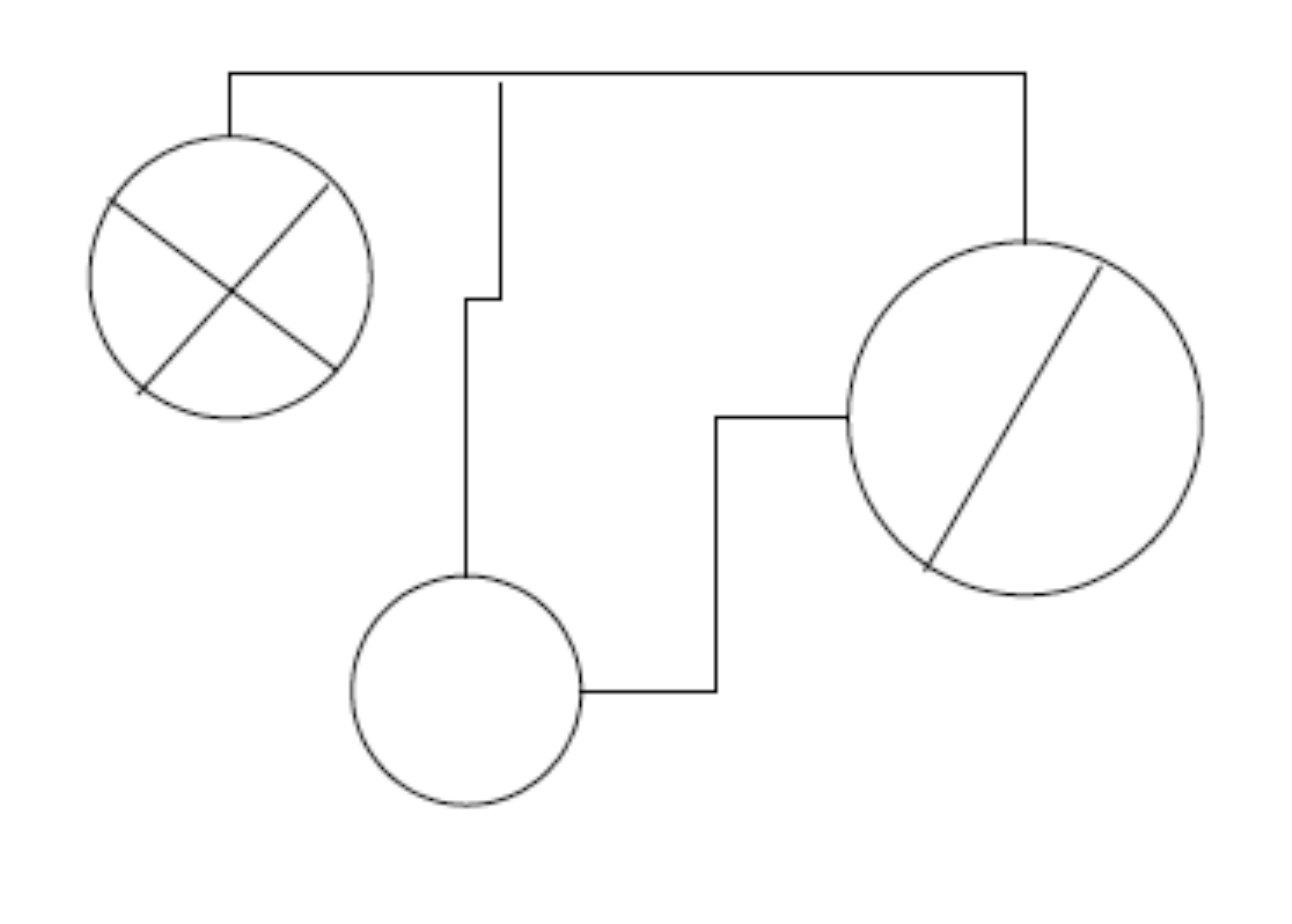}
 }

\noindent $p=3$, $r=4$ and $t =11$.

\vspace{2mm}\noindent Then we have a commutative diagram of
fibrations

$$
\xymatrix{
\mbox{Map}(S,M)\ar[d]_{\scriptstyle q_S} & \mbox{Map}(C,M) \ar[l]_{\simeq}  \ar[d]_{\scriptstyle q_C} \ar[r]& \prod_{i=1}^r
\mbox{Map}(T_i, M)\ar[d]_{\scriptstyle q_T}\\
(LM)^p\ar@{=}[r] &(LM)^p  \ar[r]^{\psi}  &M^t } $$
where $q_T$
and $\psi$ denote the evaluation maps at the extremities of the
trees and $q_C$ denotes the restriction to $\amalg_pS^1$. Since
each $T_i$ is contractible, the map $q_T$ is homotopy equivalent
to a product of diagonal maps $\Delta : M^r \to M^t$. This shows
that the fibration $q_S$ is homotopy equivalent to the pullback of
$\Delta$  along some map $\psi$.

Assume that $M$ is a 1-connected oriented Poincar\'e duality
space. Then   by Theorem A,   there exists in the derived category
of $C^*((LM)^p)$-modules  a well defined (up to homotopy) morphism
$$
(q_S)^! : C^* (\mbox{Map}(S,M))  \to C^{*-m\chi}    ((LM)^p)\,,
$$
of degree $-m\chi$ where $\chi = ( 2 -2g-p-q)$  denotes the Euler
characteristic of $S$. The {\it cohomology $(g, p+q)$-string
operation} induced by the surface $S$,
$$
H^*((LM)^q)\to H^{*-m\chi}    ((LM)^p)\,,
$$
\noindent  is
 defined as   the composition of $H^*((q_S)^!)$ with the morphism
  $H^*(\mbox{Map}(k,M))$
when $k$ denotes  the inclusion  of the outgoing circles into $S$.
For closed oriented manifold  the dual of these cohomology $(g, p+q)$-string
operations coincide with those defined by Cohen-Godin, \cite{CG}.

  In \cite{GS}
Gruher and Salvatore    prove  that the loop product   is an
oriented homotopy invariant.
 We generalize this result to all the string operations defined on a Poincar\'e
 duality space.

 \vspace{2mm}\noindent {\bf Theorem B.} {\sl If  $M$ is a 1-connected  oriented
  Poincar\'e duality space then all the $(g, p+q)$-string operations on
 $H_*(LM)$
 depend  only on the oriented homotopy type of $M$.}

\vspace{2mm}   We will now explain how   string operations can be
extended to Gorenstein spaces. For recall, a differential graded
augmented $\bk$-algebra is called  a {\it Gorenstein algebra  of dimension $d$}  if
 $  \mbox{Ext}^k_A(\bk,A)=\{0\}$ if $k\neq d$ and  $  \mbox{Ext}^d_A(\bk,A)=\{0\}$
 has dimension one.  (Here the field $\bk$ is considered as a trivial module while
 $A$ acts on itself by left multiplication.) A path connected space $M$ is
 called  a {\it  $\bk$-Gorenstein space of dimension $d$} if the nondegenerated
 singular cochain
 algebra $C^*(M)$ is a Gorenstein algebra of dimension $d$.
Gorenstein  spaces have been introduced
 in \cite{FHT} in relation with the Spivak fibration.
   Examples of Gorenstein
 spaces are given by Poincar\'e duality spaces,
 classifying spaces $BG$  of compact connected Lie groups $G$ and rational
 spaces with finite Postnikov tower. More generally,  if
  $F\to E\to B$ is a fibration in which
 $B$ is a Gorenstein space and $F$   a
 Poincar\'e duality space, then $E$ is   a Gorenstein space (\cite{FHT}).
 As an application, if $G$   acts on a path connected
 Poincar\'e duality space $M$ then   the homotopy quotient $  EG\times_GM$ is
 a Gorenstein space since it is the total space of the Borel fibration
 $M\to EG\times_GM\to BG$. String operations of those spaces have
 been recently studied using stacks and bivariant homology theory
 (\cite{BGNX}).

For   a Poincar\'e duality space $M$, the construction and the
uniqueness of the string operations on $LM$ were depending  only
on
  the existence and the uniqueness (up to homotopy), in the derived
category of $C^*(M^n)$-modules, of a map of degree $d(n-r)$,
$\Delta^! : C^*(M^r)\to C^*(M^n)$, for $r\leq n$. The main tool
to define  string operations  on  Gorenstein spaces
is the following Theorem C.

\vspace{2mm}\noindent {\bf Theorem C.}  {\sl Let $X$ be a
1-connected $\bk$-Gorenstein space of dimension $d$, and for
$r\leq n$ let $\Delta : X^r\to X^n$ be the product of diagonal
maps $X\to X^{n_i}, \sum_{i=1}^r n_i = n$. Then,
$$\mbox{Ext}_{C^*(X^n)}(C^*(X^r), C^*(X^n)) = s^{(n-r)d}\,
H^*(X^r;\bk)\,,$$ \noindent  where $C^*(X^r)$ is viewed as a
$C^*(X^n)$-module via   $\Delta $.}

\vspace{2mm} This implies that, up to homotopy and up to the
multiplication by a scalar,  there is in the derived category of
$C^*(X^n)$-modules   a unique non trivial class of morphism of
degree $(n-r)d$,

$$\Delta^! : C^*(X^r ) \to C^*(X^n)\,, $$   and for each
surface $S$ of genus $g$ with $p+q$ boundary components,
   $\Delta^!$ induces a map
$$(q_S)^! : C^*(\mbox{Map}(S,X)) \to C^*((LX)^p) \,,  $$
  in the derived category of $C^*((LX)^p)$-modules.  This map $(q_S)^!$ is uniquely
   defined up to homotopy and up to
the multiplication by a  scalar. In the same way as in the
Poincar\'e duality case, we can now defined string operations for
Gorenstein spaces. We give in section 5 some computations.   In
particular, considering the $(0, 1+2)$-string operation,  we prove:

\vspace{2mm}\noindent {\bf Theorem D.}   {\sl  Let  $\bk=\mathbb Q$ and let  $ BG$ be
the classifying space of a connected compact Lie group  then the
   loop coproduct $H_*(L(BG)) \to
H_*(L(BG)) \otimes H^*(L(BG))$ is injective.}

\vspace{2mm}  Section  1 contains brief background material on cofibrant
objects in the category of differential graded algebras and
 differential modules. In section 2, we define shriek maps for Poincar\'e duality spaces
 and prove Theorem A and B.
  In section 3   we consider string operations on   Poincar\'e duality spaces.
  Section 4 is devoted to the characteristic zero case. We make there use
  of Sullivan minimal models for making computations. In section 5 and 6
   we extend the previous definitions to Gorenstein spaces and we prove in particular Theorem C.

\section{Models for algebras and modules}

\subsection{Conventions}

All vector spaces are defined on a fixed field $\bk$ and the
unadorned $\otimes$ and $\mbox{Hom}$ mean with respect to $\bk$.
Graduations are written either as superscripts or as subscripts,
with the convention $V^k=V_{-k}$. We say that the graded vector
space $V$  has finite type if each $V^i$ is finite dimensional.
The graded dual of $V$  is denoted by $V^\#$ i.e. $V^\#_k=
(V^k)^\#$.

Differential graded algebras $(R,d)$ are assume to be   of the
form $R=\{R^k\}_{k\geq 0}$  with differential $d$ of (upper)
degree +1. If $V$ is a graded vector space, $T(V)$ denotes the
tensor algebra on $V$.

Unless we say otherwise, we  shall use the word $(R,d)$-module for
(left or right)  differential $\mathbb Z$-graded module over a
differential graded algebra $(R,d)$.  Let  $(Q,d)$ and $(Q',d)$ be
$R$-modules, we denote by $\mbox{Hom}_R(Q,Q')$ the graded vector
spaces,
$$
\mbox{Hom}^k_R(Q,Q') =\prod_{i\in \mathbb Z} \mbox{Hom}_R(Q^i ,(Q')^{i+k})
$$
with the differential $D$ defined by $D\varphi = d\circ \varphi -
(-1)^{k} \varphi \circ d$.

Recall also that the $R$-linear maps $\varphi, \psi \in
\mbox{Hom}^k_R(Q,Q')$ are {\it homotopic} if $\varphi-\psi = D
\theta $ for some $\theta \in \mbox{Hom}^{k-1}_R(Q,Q')$.

\subsection{Semifree modules in the category of $R$-modules}

\vspace{1mm}\noindent {\bf Definition.} Let $(R,d)$ a differential
graded algebra over the fixed field $\bk$. \begin{enumerate}
\item[(i)]  {\it A $R$-semifree extension of an $R$-module $(Q,d)$}  is  a
morphism of $R$-modules  of the form
$
(Q,d) \stackrel{i} \hookrightarrow (Q\otimes V, d)
$
in which
\begin{enumerate}
\item $i$ is the obvious inclusion.
\item $V= \bigcup _k V(k)$ with  $V(0)\subset \cdots  \subset V(k) \subset \cdots $ sub
vector
spaces of $V$
\item $ d(V(0)) \subset Q$ and $d(V(k) \subset Q\otimes V(k-1) $, $k\geq 1$.
\end{enumerate}
\item[(ii)] If $Q=R$ we say that {\it $(R\otimes V, d)$ is a $R$-semifree module}.
\end{enumerate}

\noindent Semifree modules are cofibrant objects in   the category
of $R$-modules. We recall here their main properties, see
\cite[Section 6]{FHT3} or \cite[§ 2]{FHT2b} for more details.

Let $ \varphi : (Q, d) \to (Q',d')$ be a  morphism of $R$-modules.

\begin{enumerate}
\item[(SF1)]  $ \varphi $ admits a   {\it relative $R$-semifree model}, i.e. there exists
a commutative diagram

\centerline{\footnotesize$
\xymatrix{
&&P\ar[d]^{\simeq}\\
Q \ar@{^(->}[urr]^{i} \ar[rr]^{\varphi} && Q'
}
$}
\noindent where $
(Q,d) \stackrel{i} \hookrightarrow (P,d) $ is a  $R$-semifree extension of $(Q,d)$.
 In particular ($Q=R$ and $\varphi=0$), every $R$-modules  admits a {\it $R$-semifree model}.

\item[(SF2)] If $P$ is a $R$-semifree module and if $\varphi$ is a
quasi-isomorphism then $\mbox{Hom}_R(P, \varphi) : \mbox{Hom}_R(P, Q)\to \mbox{Hom}_R(P, Q')$
is a quasi-isomorphism.

\item[(SF3)] Given a diagram of morphisms of $R$-modules of the form

\centerline{\footnotesize $
\xymatrix{
P\ar[d]^i\ar[r]^{\psi}&Q'' \ar[d]^{\eta}_{\simeq} \\
Q\ar[r]_{\varphi} & Q' } $} \noindent where $i:Q\hookrightarrow P$
is a $R$-semifree extension and $\eta$ is a quasi-isomorphism,
then there exists a morphism of $R$-modules $\psi' : Q \to Q''$
(unique up to homotopy) such that $\psi' i=\psi$ and $\eta\psi'
\sim \varphi$.

\end{enumerate}

\noindent By definition, if $P$ is a semifree model of $N$, then,
 $$
\mbox{Ext}^n_R(N,Q):= H^n\left(  \mbox{Hom}_R( P, Q)\right)$$
while if  $ (Q',d) $ is a right  module,
  $$ \mbox{Tor}_n^R(Q',N):= H_n\left(   Q' \otimes_R  P\right)\,.
$$

\subsection{Bar constructions and Hochschild chain complex}

 Let $A=\bk \oplus \bar A $ be a supplemented graded differential graded
algebra, $N$ a right $A$-module and $M$ a left $A$-module. The
two-sided (graded) bar construction on $A, M$ and $N$, $\mathbb B
(N, A, M)$, is the differential bigraded vector space $ \mathbb B
(N, A, M) =  (N \otimes T (s \bar A ) \otimes M, d)
 $ defined as follows:

 For $k\geq 1$, we write  $n[a_1|a_2|...|a_k]m =
 n\otimes sa_1\otimes \ldots \otimes sa_k\otimes m \in
\mathbb B_k (N, A, M)
 $.  If $k=0$, we write $n[\,]m = n\otimes 1\otimes m \in
N \otimes T^0(s \bar A)\otimes M$. The
 differential $d= d_0+d_1$ is defined by
 $$
 \renewcommand{\arraystretch}{1.6}
 \begin{array}{rll}
 d_0  ( n[a_1|a_2|...|a_k]m)& = d( n) [a_1|a_2|...|a_k]m  -
\displaystyle\sum
 _{i=1}^k (-1)^{\epsilon _i}  n[a_1|a_2|...|d(
 a_i)|...|a_k]m\\ &+ (-1) ^{\epsilon _{k+1}}  n[a_1|a_2|...|a_k]d(m)
 \\[3mm]
 d_1  (n[a_1|a_2|...|a_k]m) &= (-1) ^{| n|}  na_1[a_2|...|a_k]m +
 \displaystyle\sum _{i=2}^k (-1) ^{\epsilon _i}
 n[a_1|a_2|...|a_{i-1}a_i|... | a_k]m \\
 &- (-1)^{\epsilon _{k}}  n[a_1|a_2|...|a_{k-1}]
 a_k m
 \end{array}
 \renewcommand{\arraystretch}{1}
 $$
  Here
 $\epsilon _i = | n| + \sum _{j<i}  |s a_j| $.

The (graded) bar construction on $A$ is the differential graded coalgebra
  $
\overline{\mathbb B} A =\mathbb B(\bk ,A,\bk )
 $,  (see \cite[(4.2)]{FHT2b}),   whose comultiplication $\phi$ is defined by
 $\phi ([a_1\vert \cdots \vert a_r] )=
 \sum_{i=0}^r \, [a_1\vert \cdots \vert a_i]\otimes [a_{i+1}\vert \cdots   \vert a_r]$.

 For each $A$-module $M$, $\mathbb B(\bk, A,M)$ is   a
$\overline{\mathbb
 B}(A)$-comodule and  the morphism $\varphi :
\mathbb B (A,A,M) \to M$
 defined by $\varphi(a[\,]m) = am$
 is a $A$-semifree
 model for $M$, \cite[Lemma 4.3]{FHT2b}.

 Denote by $A^{op}$ the opposite algebra and by $A^e = A \otimes
A^{op}$.  The bar construction $\mathbb B (A,A,A)$ is a $A^e$-semifree
model of $A$.  The {\it Hochschild chain complex of $A$} is
the complex $C\!H_*(A) = (A,d)\otimes_{A^e}(\mathbb B(A,A,A),d)$
where $A$ is viewed as a $A^e$-module via the multiplication. Its
homology $H\!H_*(A)= \mbox{Tor}_{A^e}(A,A) $ is called the {\it Hochschild homology of $A$
(with coefficients in $A$)}.

\subsection{Free models}

 \noindent {\bf Definition.}
\begin{enumerate}
\item[(i)]  A  {\it free extension of a differential graded algebra $(A,d)$}  is  a
morphism of differential graded algebras   of the form $ (A,d)
\stackrel{i} \hookrightarrow (A\coprod  T(V), d) $ in which
\begin{enumerate}
\item $i$ is the obvious inclusion.
\item  $V= \bigcup _k V(k)$ with  $V(0)\subset \cdots  \subset V(k) \subset \cdots $ subspaces of $V$
\item $ d(V(0)) \subset A$ and $d(V(k) \subset A\coprod T( V(k-1)) $, $k\geq 1$.
\end{enumerate}
\item[(ii)] If $A=\bk$ we say that {\it $(R,d)$ is a  free extension}.
\end{enumerate}

\noindent Free extensions  are cofibrant objects in  the category
of differential graded algebras.  We recall here  some of their
main properties, see \cite{HL}, \cite{FHT2b} or \cite{FHT2} for
more details. Recall only  that    two  morphisms of differential
graded algebras $ \varphi, \psi : T(V) \to R$  are {\it homotopic}
if there exists a $(\varphi,\psi)$-derivation $\theta$ such that
$\varphi-\psi = d\circ \theta -\theta \circ d$.

Let $(A,d)$ and let $(A',d)$ be differential graded algebras  and
let $ \varphi : (A, d) \to (A',d')$ be  a homomorphism of
differential graded algebras .

\begin{enumerate}
\item[(FE1)]   $ \varphi $  admits a {\it relative   free model}, i.e.
there exists a commutative diagram

\centerline{\footnotesize$
\xymatrix{
&& (A \coprod T(V),d) \ar[d]^{\simeq}\\
(A,d) \ar@{^(->}[urr]^{i} \ar[rr]^{\varphi} && (A',d) } $}
\noindent where $ (A,d) \stackrel{i} \hookrightarrow   (A \coprod
T(V)$ is a free extension of $(A,d)$.
 In particular   every  differential graded algebra  admits a  free model.

\item[(FE2)]  Given a homotopy commutative  diagram  in the category of differential graded
algebras of the form

\centerline{\footnotesize$
\xymatrix{
(T(V) ,d)\ar[d]^i\ar[rr]^{\psi} && (A,d)\ar[d]^{\eta}_{\simeq}\\
(T(V\oplus W), d )   \ar[rr]^{\varphi}&& (A',d) } $} \noindent
where $i$  is a free extension of $(T(V),d)$ and $\eta$ a
quasi-isomorphism,  then there exists a  morphism of differential
graded algebras (unique up to homotopy)   $\psi': (T(V\oplus W), d
) \to (A,d)$ such that $\psi'i=\psi$ and $\eta\psi'\sim \varphi$.

Moreover, if $\eta$ is surjective and the diagram commutative, we
can choose $\psi'$ such that $\eta\psi'= \varphi$.

\item[(FE3)] If $\eta : (T(V),d) \to A$ and $\eta' : (  T(W),d) \to A'$
are free models of $(A,d)$, then   there is a quasi-isomorphism
$\alpha : (T( W),d)  \to (T(V),d) $ such that $\eta \circ \alpha$
is homotopic to   $\eta'$.

\end{enumerate}


\section{Poincar\'e duality spaces}

\subsection{Shriek maps}

 Let   $M$ be an oriented Poincar\'e duality space of dimension $m$,
  $X$ any path connected space with homology of finite type  and  $f : X \to
M$ be a  map, then $C^*(X)$ is a $C^*(M)$-module via $f$.
Therefore, for each  integer  $q$ we have a canonical  linear map
$$
\theta : \mbox{Ext}_{C^*(M)}^{q}(C^*(X), C^*(M)) \to
\mbox{Hom}^q(H^{m-q}(X), H^m(M)) \,, \varphi \mapsto [\varphi]\,.
$$

\vspace{2mm}\noindent {\bf Lemma 1.} (\cite{LS2})  {\sl With the
above notations, the linear map   $\theta $ is an isomorphism. }

\vspace{2mm}\noindent {\bf Proof.} Let $\psi : R \to C^*(X)$ be a
$C^*(M)$-semifree model of $C^*(X)$. Then  by property (SF-2) of
semifree models, the cap product with a cycle representing the
orientation class of $M$, $ cap_M : C^*(M) \to C_*(M)$, induces a
quasi-isomorphism
$$
 \mbox{Hom}_{C^*(M)} (R  , C^*(M)) \to
\mbox{Hom}_{C^*(M)}( R, C_*(M))\,.
$$

Now let us consider a linear basis, $(\eta_i)_i$, of
$\mbox{Hom}(H^{m-q}(X), H^m(M))$ and the corresponding basis
$([a_i])_i$ of $H_{m-q}(X)$ under the  linear isomorphism
$$
 \mbox{Hom}(H^{m-q}(X), H^m(M) \cong  \mbox{Hom}(H^{m-q}(X), \bk)\cong H_{m-q}(X)\,.
 $$
By definition $a_i$ is a cycle in $C_{m-q}(X)$ and  we denote by
$cap_{a_i} :C^\ast(X) \to C_\ast(X)$ the cap product by $a_i$. We
deduce, for each $i$,  a homotopy commutative diagram in the
category of $C^*(M)$-modules of the form
$$
\xymatrix{
R\ar[d]_{\psi}\ar[rr]^{\psi_i}&&C^*(N)\ar[d]^{cap_M}_{\simeq}&\ar@{}[d]_{.}\\
C^*(X) \ar[r]^{cap_{a_i}}&C_*(X)\ar[r]^{C_*(f)}&C_*(M)& }
$$
Since $ {\theta} (\psi_i)=\eta_i$, $ {\theta}$ is surjective.  In
fact this is a bijection since the two vector spaces have the same
dimension. This last fact comes from the Cartan-Eilenberg
associativity formulae. Indeed,
$$
\begin{array}{ll}
\mbox{Ext}^{q}_{C^*(M)}(C^*(X), C^*(M))&\cong \mbox{Hom}
\left(\mbox{Tor}_{q}^{C^*(M)}(C^*(X),
 C_*(M)),\, \bk\right)\\
 &\cong
\mbox{Hom}\left(\mbox{Tor}_{q-m}^{C^*(M)}(C^*(X), C^*(M)),\bk\right)\\
&\cong \mbox{Hom}(H^{m-q}(X),\bk)\,.
\end{array}
$$
\rightline{$\square$}

\vspace{2mm}  Lemma 1 implies directly the next result.

\vspace{2mm}\noindent {\bf Theorem 1.} {\sl  Let $M$ and $N$ be
two Poincar\'e duality spaces and  $ f : N \to M$ be a map. Then
there is,  in the derived category of $C^*(M)$-modules,  a unique
map (up to homotopy) $$f^{!} : C^*(N) \to C^*(M)$$ such that
$H^*(f^!)$ maps the fundamental class of $N$  to the fundamental
class of $M$.}

\vspace{2mm} Now let consider diagram (*) and hypothesis (H) from
the Introduction. Then,

\vspace{2mm} \noindent {\bf Theorem 2.} {\sl  In the derived
category of $C^*(E)$-modules there exists a unique map (up to
homotopy) $g^{!} : C^*(E')\to C^*(E)$ making commutative the
diagram
$$\xymatrix{
C^*(E') \ar[rr]^{g^{!}}&& C^*(E)&\ar@{}^{.}\\
C^*(N)\ar[u]^{\scriptstyle C^*(p')}\ar[rr]^{f^!}&&
C^*(M)\ar[u]_{\scriptstyle C^*(p)}& }$$   }

\vspace{2mm}\noindent {\bf Proof.} Denote by  $(P,d)\to C^*(N)$ a
$C^*(M)$-semifree model of $C^*(N)$  and consider $C^*(E)$ as a
right $C^*(M)$-module via $C^*(p)$. Since $M$ is simply connected,
by (\cite{FHT3}, Theorem 7.5), there exists a quasi-isomorphism of
$C^*(E)$-modules

$$C^*(E)\otimes_{C^*(M)}(P,d) \to C^*(E')
$$ that makes $C^*(E)\otimes_{C^*(M)}(P,d) $ a $C^*(E)$-semifree
model of $C^*(E')$.

Now let $f^! : (P,d) \to C^*(M)$ be a shriek map associated to $f$
by Theorem 1. The tensor product
$$
g^! = 1 \otimes f^! : C^*(E) \otimes_{C^*(M)}(P,d) \to
C^*(E)\otimes_{C^*(M)}C^*(M)
$$
 is a morphism of left $C^*(E)$-modules making commutative the diagram
$$
\xymatrix{ C^*(E) \otimes_{C^*(M)}(P,d)  \ar[r]^{g^!}&
C^*(E)\otimes_{C^*(M)}C^*(M)\ar[r]^(0.7){\cong} &C^*(E)\\
(P,d)\ar[u]^{\scriptstyle h} \ar[rr]^{f^!} &&
C^*(M)\ar[u]_{\scriptstyle C^*(p)} } $$
 with $h(a) = 1
\otimes a$, $a\in R$.

We write  $(P,d) = (C^*(M)\otimes V,d)$ with $d(V(q))\subset
C^*(M)\otimes V(q-1)$. Then $C^*(E)\otimes_{C^*(M)}(P,d)\cong (
C^*(E)\otimes V,d)$ with $d(V(q))\subset C^*(E)\otimes V(q-1)$.
Since $g^!$ is completely determined by its values on $V$, there
is only one way to extend $f^!$ into a map $g^!$ of
$C^*(E)$-modules. \hfill $\square$

\vspace{2mm} Consider once again the diagram (*) with the
hypothesis (H). Denote by $(R,d)$   a semifree model for $C^*(N)$.
We can then associate to each morphism of $C^*(M)$-module $k :
R\to C^*(M)$ the morphism of $C^*(E)$-module $1\otimes k :
C^*E\otimes_{C^*M}R \to C^*E$. This correspondence  defines a
linear map
$$\Phi : \mbox{Ext}_{C^*(M)}(C^*(N), C^*(M)) \to
\mbox{Ext}_{C^*(E)}(C^*(E'), C^*(E)))\,.$$ Clearly from the
definition of $g^!$, we have
$$g^!=\Phi (f^!)\,.$$

\subsection{Homotopy invariance of the shriek maps}

Consider now the following commutative diagram
$$
\xymatrix{
  &&&
X  \ar[rr]^{g} \ar'[d][dd]_{q}\ar[dl]^{k'} && E \ar[dd]^{p}
\ar[dl]^{k}
\\
 && X' \ar[rr]^(0.6){g' } \ar[dd]_{q'} && E'
\ar[dd]^(0.3){p'}&& \\
 &&&N \ar'[r][rr]_(-0.3){f} \ar[dl]^{h'}
&& M\ar[dl]^{h} \\
 &&N'\ar[rr]_{f'} &&M'
 }
$$
where
$$\left\{ \begin{array}{l}
\bullet \hspace{5mm}  \mbox{$M$ and $M'$ are 1-connected}\\
\bullet \hspace{5mm} \mbox{$M$, $M'$, $N$ and $N'$ are Poincar\'e
duality complexes}\\
\bullet \hspace{5mm}\mbox{$H^*(E)$ and $H^*(E')$ are finite type
vector spaces}\\
\bullet \hspace{5mm}\mbox{The vertical maps are fibrations}\\
\bullet\hspace{5mm}\mbox{$q$ and $q'$ are respectively pullback of
$p$ and $p'$ along $f$ and $f'$}\end{array}\right.$$

 \vspace{2mm}\noindent {\bf Theorem 3.} {\sl With the above notations, if
$h, h', k$ and $k'$ are homotopy equivalences, then,
$$H^*(g^!)\circ H^*(k') = H^*(k)\circ H^*((g')^!)\,.$$}

\vspace{2mm}\noindent {\bf Proof.}  The naturality of the
isomorphism $\theta$ described at the beginning of section 2.1
gives the following commutative diagram of isomorphisms
$$
\xymatrix{ \mbox{Ext}_{C^*M}(C^*N, C^*M) \ar[r]^{\rho}
\ar[d]^{\scriptstyle \theta}  &
\mbox{Ext}_{C^*M'}(C^*N', C^*M) \ar[d]^{\scriptstyle \theta}  & \mbox{Ext}_{C^*M'}(C^*N', C^*M') \ar[d]^{\scriptstyle \theta}  \ar[l]_{\sigma}&\\
\mbox{Hom}(H^{m-n}N, H^mM) \ar[r]^{ \overline{\rho}}  & \mbox{Hom}
(H^{m-n}N', H^mM)& \mbox{Hom}(H^{m-n}N',
H^mM')\ar[l]_{\overline{\sigma}}& }$$ Here $\rho$ and $\sigma$ are
the evident isomorphisms,   $\overline{\rho}\theta
(f^!)=H^*(f^!)\circ H^*(h')$ and $\overline{\sigma}\theta
(f'!)=H^*(h)\circ H^*(f'!)$. Now $\sigma (f'!)=\rho (f^!)$ because
their images by $\theta$ coincide.   Now by the naturality of
$\Phi$ we have a commutative diagram of graded vector spaces where
$\rho'$ and $\sigma'$ are the evident isomorphisms induced by the
homotopy equivalences $k$ and $k'$.
$$
\xymatrix{ \mbox{Ext}_{C^*E}(C^*X, C^*E) \ar[r]^{\rho'}&
\mbox{Ext}_{C^*E'}(C^*X', C^*E) &
\mbox{Ext}_{C^*E'}(C^*X', C^*E') \ar[l]_{\sigma'}&\ar@{}[d]^{.}\\
\mbox{Ext}_{C^*(M)}(C^*N, C^*M) \ar[r]{ \rho} \ar[u]^{\scriptstyle
\Phi} & \mbox{Ext}_{C^*(M')}(C^*N',
C^*M){\sigma}\ar[u]^{\scriptstyle \Phi} &
\mbox{Ext}_{C^*(M')}(C^*N',
C^*M')\ar[l]_{\sigma}\ar[u]^{\scriptstyle \Phi}& }$$  In
particular $\rho'(g^!)= \sigma'({g'}^!)$. This implies the
result.\hfill $\square$

\subsection{Naturality  of the shriek maps}

\vspace{3mm} We consider a diagram of fibrations over $B$,
$$\xymatrix{
 X \ar[rr]^{g} \ar[dr]_{\scriptstyle p}&& Y\ar[dl]^{\scriptstyle q}\\
&B\,.}
$$
 We form the pullback
of this triangle  along  a map $f : B'\to B$:
$$
\xymatrix{
  &&&
X' \ar[rr]^{h} \ar'[d][dd]^{p'}\ar[dl]^{g'} && X \ar[dd]^{p}
\ar[dl]^{g}&\ar@{}[dd]^{.}
\\
 && Y' \ar[rr]^(0.6){k } \ar[dr]_{q'} && Y
\ar[dr]^(0.3){q}&& \\
 &&&B' \ar'[rr]^{f}
&& B &
 }
$$

\vspace{2mm}\noindent {\bf Theorem 4.}  {\sl  If $B$ and $B'$ are
 1-connected oriented Poincar\'e duality spaces, and $H^*(X)$ and
$H^*(Y)$ of finite type,  then,
$$H^*(g)\circ H^*(k^!)  = H^*(h^!) \circ H^*(g'): H^*(Y')\to H^*(X)\,.$$}

\vspace{2mm}\noindent {\bf Proof.} We put  $(R,d) = C^*(B)$ and we
denote by $\varphi : (R,d) \to (R\otimes V,D)$ an $R$-semifree
model of $C^*(f)$.  We denote also by $\psi : (R_1,d) \to (R_2,d)$
a model of $g$ in the category of  differential graded
$R$-algebras. Then a model of $g'$ in the category of
$R_1$-modules is
$$\psi\otimes 1 : (R_1,d)\otimes_R(R\otimes V,D) \to
(R_2,d)\otimes_R(R\otimes V,D)\,.
$$ Now denote by $\theta : (R\otimes V,D)\to (R,d)$  a
representative for $ f^!$. The maps $C^*(g)\circ k^!$ and $
 h^!\circ C^*(g')$ are given respectively by
 $$(R_1,d)\otimes_R(R\otimes V,D) \stackrel{1\otimes\theta}{\to} (R_1,d)\otimes_R (R,d) = (R_1,d)\stackrel{\psi}{\to}
 (R_2,d)\,.$$
 and
 $$(R_1,d)\otimes_R(R\otimes V,D) \stackrel{\psi\otimes 1}{\to} (R_2,d)\otimes_R(R\otimes V,D)
 \stackrel{1\otimes \theta}{\to} (R_2,d)\,.
 $$
Since   $C^*(g)\circ k^!$ and $
 h^!\circ C^*(g')$ are homotopic, they  induce  the same map in
 cohomology. \hfill $\square$

\subsection{Integration along the fiber}

Let $F \to E \stackrel{p}{\to} B$ be a  fibration in which $F$ is
an oriented  Poincar\'e duality space of cohomological dimension
$d$ and $B$ a 1-connected space. The integration along the fiber
is the linear map, $\int_F: H^*(E) \to H^{*-d}(B)$,   defined  as
the composition of natural maps arising in the Serre spectral
sequence of the fibration

\centerline{\footnotesize $
 H^*(E) \twoheadrightarrow E^{*,d}_\infty \subset E_2^{*,d} = H^*(B)\otimes
H^d(F) \cong H^{*-d}(B)\,.$}

\vspace{2mm}\noindent {\bf Theorem 5.} {\sl With the notation
above, if $E$ and $B$ are oriented Poincar\'e duality spaces, then
$H^*(p^!) : H^*(E) \to H^{*-d}(B) $   is the integration along the
fiber.}

\vspace{2mm}\noindent {\bf Proof.} The differential graded algebra
$C^*(E)$ admits a semifree $C^*(B)$-module of the form
$(C^*(B)\otimes H^*(F),d)$ with $d(H^r(F) ) \subset C^*(B)\otimes
H^{<r}(F)$ (\cite{FHT}, Lemma A.3). In this setting  it is clear
that the integration along the fiber is the map induced in
homology by the composite

\centerline{\footnotesize $ \xymatrix{
 C^*(B)\otimes H^*(F)\ar@/^2pc/[rr]^{\varphi} \ar@{>>}[r] & C^*(B) \otimes
 \left( H^*(F)/H^{<d}(F)\right)   \ar[r]^(0.7){\cong} &
 C^{*-d}(B)
 }
 \,.$}

\vspace{1mm}\noindent
  We remark   that $\varphi$ is a morphism of
$C^*(B)$-modules of degree $-d = b-e$ where $b=$ dim$\,B$ and $e=$
dim$\,E$ . Moreover the image by $H^*(\varphi)$ of    the
fundamental class of $E$ is the fundamental class of $B$. The
results follows now from the uniqueness property in Theorem 1.

\hfill $\square$

\vspace{2mm}\noindent {\bf Theorem 6.} {\sl If the fibration $F
\to E \stackrel{p'}{\to} B$ is the pullback, along a map $f : B
\to M$, of a fibration $F \to X\stackrel{p}{\to} M$ in which all
the involved spaces are Poincar\'e duality spaces, then
$H^*((p')^!)$ is the integration along the fiber.}

\vspace{2mm}\noindent {\bf Proof.} Let $(C^*(M)\otimes H^*(F),d)$
be a semifree model for $C^*(X)$. Then the tensor product
$(C^*(B)\otimes_{C^*(M)}\left( C^*(M)\otimes H^*(F)\right),d)$ is
a semifree model for $C^*(E)$ as a $C^*(B)$-module. The respective
integrations along the fiber are then given by the morphisms
$\varphi$ and $\varphi'$ in the following commutative diagram,
$$
\xymatrix{
 C^*(M)\otimes H^*(F)\ar@/^2pc/[rr]^{\varphi} \ar@{>>}[r] \ar[d]_{\scriptstyle
 C^*(f)\otimes 1}& C^*(M) \otimes \left( H^*(F)/H^{<d}(F)\right)
 \ar[r]^(0.7){\cong} &  C^{*-d}(M) \ar[d]^{\scriptstyle C^*(f)} &\ar@{}[d]^{.} \\
  C^*(B)\otimes H^*(F)\ar@/_2pc/[rr]^{\varphi'}\ar@{>>}[r]  &C^{*}(B)
\otimes \left( H^*(F)/H^{<d}(F)\right) \ar[r]^(0.7){\cong}&
C^{*-d}(B) &} $$ Since $\varphi$ and $\varphi'$ are respectively
morphisms of   $C^*(M)$ and $C^*(B)$-modules of degree dim$\,M-$
dim$\,X$,   the uniqueness property of Theorem 2  implies that
$\varphi'$ is a representative of the shriek map $(p')^!$. \hfill
$\square$

\subsection{The  intersection map  with the fiber}

Consider the diagram $(\dagger)$ of the introduction.
$$
\xymatrix{
F \ar[r]^{j} \ar[d] & E\ar[d]&\ar@{}[d]^{.}\\
\{b_0\} \ar@{^(->}[r]^{i} & M& }
$$
By Theorem 2, in the derived category of $C^*(E)$-modules there
exists  a unique morphism (up to homotopy)   $j^!$ making
commutative the diagram,

 \centerline{\footnotesize$ \xymatrix{
C^*(F) \ar[r]^{j^!}& C^*(E)&\ar@{}[d]^{,}\\
\bk \ar[u] \ar[r]^{i^!}& C^*(M)\ar[u]&
}$}

 \noindent where $H^*(i^!)$ is   the multiplication by the fundamental class
of $M$.  The morphism $H^*(j^!)$, is called   the intersection map
with the fiber,  and the purpose of this section is the
description of a model of $H^*(j^!)$ in terms of cochain
complexes.

 Let $(T(V),d)$ be a  free model of $C^\ast(M)$ ( cf section 1.4)  and
 $i :  (T(V),d)\to  (W\otimes T(V),d)$ be a right
$(T(V),d)$-semifree model for the   fibration $p : E\to M$, i.e.,
there is a commutative diagram of right $T(V)$-modules, where the
vertical lines are quasi-isomorphisms
$$
\xymatrix{  (T(V),d)  \ar[rr]^{\rho} \ar[d]_{\scriptstyle \simeq} &&
 (P\otimes T(V),d)\ar[d]^{\scriptstyle \simeq}\\
C^*(M) \ar[rr]^{C^*(p)}  && C^*(E) \,. }$$ Since $M$ is
1-connected, we may suppose that $V=V^{\geq 1}$. We denote   by $q
: (T(V),d)  \twoheadrightarrow (A,d)$ a surjective
quasi-isomorphism of differential graded algebras, with
    $A^0=\bk$ and $A^{>m}=0$. Using $q$ we  define the cochain complexes

$$(P, \overline{d}): = (P\otimes T(V),d)\otimes_{T(V)}
 (\bk,0)\hspace{3mm}\mbox{and}\, (P\otimes A,d):=( P\otimes
T(V),d)\otimes_{T(V)}(A,d)\,.$$    Finally  we denote by $\omega$
a cocycle in $A$ representing the fundamental class of
$H(A)=H^\ast(M)$.

\vspace{2mm}\noindent {\bf Theorem 7.} {\sl  With the notation
above, the map of  right  $A$-modules
$$  (P, \overline{d})  \to ( P
 \otimes A,d)\,, \hspace{1cm} x\mapsto x\otimes\omega\,,$$
 \noindent   is a
representative of $j^!$. }

\vspace{2mm}\noindent {\bf Proof.}  Let $k : (T(V),d)
\hookrightarrow (T(V\oplus W),d)$ be a  relative free model for
$C^*(p): C^*(M)  \to   C^*(E)$, and  $(T(V),d) \hookrightarrow
(T(V\oplus E),d)$ be a  relative free model for $ C^*(i): (C^*(M)
 \to  C^*(\{b_0\}) $.
 $$
\xymatrix{  (T(V),d) \ar[rr]^{k} \ar[d]_{\scriptstyle \simeq}&&
(T(V\oplus W),d) \ar[d]^{\scriptstyle \simeq}&  (T(V),d) \ar[rr]
\ar[d]_{\scriptstyle \simeq}&& (T(V\oplus E),d)
\ar[d]^{\scriptstyle \simeq}\\
C^*(M) \ar[rr]^{C^*(p)}  && C^*(E)\,, &  C^*(M) \ar[rr]^{C^*(i)}
&& C^*(\{x\}) \,.} $$
 Since $(T(V\oplus E),d)$ is a
left $(T(V),d)$-semifree module, $i^!\in
\mbox{Ext}_{C^*(M)}^d(\bk, C^*(M))\cong \mbox{Ext}_{T(V)}^d(\bk,
T(V))$ can be represented by a morphism of $T(V)$-modules of
degree $q$,
$$\nu : (T(V\oplus E),d) \to (T(V),d)\,.$$
Therefore by Theorem 2, $j^!$ is represented by the morphism of
left $(T(V\oplus W),d)$-modules

\centerline{\footnotesize $1\otimes \nu : (T(V\oplus W),d)
\otimes_{T(V)}(T(V\oplus E),d) \to
 (T(V\oplus W),d)\otimes_{T(V)}(T(V ),d) = T(V\oplus W),d)\,.
 $}
\noindent Now observe that the projection $q : (T(V),d) \to (A,d)$ makes
$(A,d)$ a $(T(V),d)$-module and we have a commutative diagram of
$(T(V),d)$-modules

\centerline{\footnotesize$
\xymatrix{
(T(V\oplus E),d) \ar[rr]^{\nu} \ar[d]_{\scriptstyle \varepsilon}&&
(T(V),d)\ar@{->>}[d]_{\simeq}^{\scriptstyle
q}& \\
(\bk, 0) \ar[rr]^{\varphi}&& (A,d)\,.&}$} \noindent  where
$\varphi (1) = q\nu (1)$.

Since $(T(V\oplus W),d)$ is a $(T(V),d)$-semifree module, taking
the tensor product with $(T(V\oplus W),d)$ over $T(V)$ yields  a
commutative diagram of $(T(V\oplus W),d)$-modules where the
vertical lines are quasi-isomorphisms

\centerline{\footnotesize$
\xymatrix{
(T(V\oplus W),d)\otimes_{T(V)} (T(V\oplus E),d)
\ar[d]^{\simeq}_ {\scriptstyle 1\otimes \varepsilon} \ar[rr]^{1\otimes \nu}&
& (T(V\oplus W),d)\ar[d]^ {\scriptstyle 1\otimes q}_{\simeq} & \\
 (T(V\oplus W),d)\otimes_{T(V)} \bk
 \ar[rr]^{1\otimes\varphi}&& (T(V\oplus W),d)\otimes_{T(V)}
 (A,d)\,,& }$}
 \noindent
 and $1\otimes \varphi$ is a new representative for $j^!$.

 Since two semifree models of a differential module are always
 quasi-isomorphic, we have a
 quasi-isomorphism of $(T(V),d)$-semifree modules
$$\theta : (P\otimes T(V) ,d) \to T(V\oplus W),d)\,.$$
  Thus tensoring by  $(A,d)$  yields a quasi-isomorphism of
 $(A,d)$-semifree modules
$$\theta\otimes 1: ( P\otimes A,d) \to
 (T(V\oplus W),d)\otimes_{T(V)}(A,d)\,.$$
 \noindent Now the commutativity of the diagram

 \centerline{\footnotesize$
 \xymatrix{
(T(V\oplus
 W),d)\otimes_{T(V)}\bk \ar[rr] ^{1\otimes \varphi}&&
(T(V\oplus W),d)\otimes_{T(V)} (A,d) \\
 (P,\overline{d}):= (P\otimes A,d)\otimes_A \bk  \ar[u]^{\scriptstyle \theta\otimes 1 }
 \ar[rr]^{1\otimes\varphi} && (P\otimes A,d)\ar[u]_{\scriptstyle \theta\otimes 1 } }$}
 \noindent shows that the map $1\otimes \varphi  : (P,\overline{d} ) \to (P\otimes
 A,d)$
 is a representative for $j^!$.

 Now note that $\omega$ and
 $\varphi(1)$ represent the fundamental class, the map  $\varphi$ is
 thus homotopic to the map $\varphi': (\bk, 0)\to (A,d)$ defined
 by
$ \varphi'(1)=\omega$. The maps  $1\otimes \varphi$ and $1\otimes
\varphi'$ are then also homotopic and both   represent
$j^!$.\hfill $\square$

\subsection{The integration map with the fiber in the  free loop space fibration}

  Let   $M$ be a simply connected oriented Poincar\'e
duality space, and   $\Omega M\to LM\to M$ be the free loop space
fibration. We denote by $(A,d)$ a differential graded algebra
quasi-isomorphic to $C^*(M)$ and satisfying $A^0=k$ and
$A^{>m}=0$, with $m=$ dim$\,M$.

\vspace{2mm}\noindent {\bf Theorem 8.}  {\sl If  $\omega\in A^m$
is a representative of the fundamental class of $M$, then there is
a commutative diagram

\centerline{ \footnotesize $
\xymatrix{
H^*(\Omega M)\ar[d]^{\cong}  \ar[rr]^{H(j^!)} &&H^\ast(LM)\ar[d]^{\cong} \\
H(\overline{\mathbb B}(A))\ar[rr]&&  H\!H_*(A) } $}
 \noindent where the lower map is induced by the multiplication by
 $\omega$,
$$ \overline{\mathbb B}
(A) \to C\!H_*(A) \,, \quad [a_1\vert \cdots \vert a_k] \mapsto
\omega[a_1\vert \cdots \vert a_k]\,.$$ }

\vspace{2mm}\noindent {\bf Proof.} Let $K$ be a 1-reduced finite
 simplicial set of dimension $m$ homotopy equivalent to  $M$. We denote by
$C_*K $   the normalized chain complex on $K$ and by $\omega\in
C^m(K)$  a representative of the fundamental class. Clearly we
have only to prove the result for $(A,d) = C^*(K)$. In \cite{S} Szczarba
gives an explicit quasi-isomorphism of chain
algebras $\alpha  : \Omega C_*K \to C_*(GK)$ where $\Omega$ denote
the cobar construction and $GK$ the loop group associated to $K$.
Since $C_*(K)$ is of finite type, by duality we get a
quasi-isomorphism of coalgebras $\alpha^\#: C^*(GK)\to
\overline{\mathbb B}(C^*K)$. By dualizing the results of Hess,
Parent and Scott (\cite{HPS}, Theorem 4.4 and Theorem 5.1), there
is a   multiplication (associative only up to homotopy) on the
Hochschild   complex $C\!H_*(C^*K)$ and a commutative diagram,

\centerline{ \footnotesize $
\xymatrix{
C^*(GK) \ar[rr]^{\alpha^\#}&& \overline{\mathbb
B}C^*(K)&\\
C^*(LK)\ar[u]^{\scriptstyle q} \ar[rr]^{\tau}&& C\!H_*(C^*K)\ar[u]&,\\
C^*K \ar@{=}[rr] \ar[u]^{\scriptstyle C^*(p)} &&C^*K \ar[u]&}
$}
\noindent where $\tau$ is a quasi-isomorphism
preserving  the products on $C^*(LK)$ and $CH_*(C^*K)$, up
to homotopy. Moreover $\tau (\,\omega\cup b\,) = \omega\otimes
\alpha^\# q(b)$ for $b\in C^+(LK)$. Denote now by $\varphi :
(P\otimes  C^* K,d)  \to C^*(LK)$  a semifree resolution of
$C^*(LK)$ as $C^*K$-module. Then the composition $\tau\circ
\varphi : (P\otimes C^*K,d)\to C\!H_*(C^*K)$ makes commutative the
diagram

\centerline{ \footnotesize $
\xymatrix{
(P,\overline{d})\ar[d]_{\scriptstyle \psi} \ar[rr]^{\simeq}&&\overline{\mathbb
B}(C^*K)\ar[d]^{\scriptstyle \psi'}&\ar@{}[d]^{,}\\
 (P\otimes C^*K,d)\ar[rr] ^{\tau\circ \varphi}&&
 C\!H_*(C^*K)&
 }$}

 \vspace{1mm}\noindent where $\psi$ and $\psi'$ consist in the left multiplication by
 $\omega$.
 Since the horizontal arrows in this diagram   are quasi-isomorphisms, we deduce
  that the multiplication by $\omega$, $\psi' :\overline{\mathbb
B}(C^*K)\to C\!H_*(C^*K)$ is a model for $j^!$ and induces in
cohomology the  intersection map with the fiber $H^*(j^!)$. \hfill
$\square$

\section{String topology}

 \subsection{Homotopy invariance of string operations}

 Since, as explained in the introduction the a $(g, p+q)$string operation in
 cohomology  is defined by composition of $H^*(q^!_S)$ with $H^*(c)$ where $c$ is
 a natural map. It follows directly from Theorem 3 that  a $(g, p+q)$string operation
  in cohomology is invariant with respect to orientation preserving homotopy
  equivalences between Poincar\'e duality spaces. Thus Theorem B of the Introduction is proved.

\subsection{The loop product}

The loop product, or  $(0,2+1)$-string operation, can be described
as follows. The injection of $S^1 \amalg S^1$ into the "pair of
pants" surface is homotopy equivalent to the injection of $S^1
\amalg S^1$ into the subspace $X \subset \mathbb R^2$ that is the
union of the interval $[0,1]\times\{0\}$ with the two circles of
radius $1/4$ centered at the points $(0, 1/4)$ and $(1, 1/4)$.
Since this injection is a cofibration, for any 1-connected
oriented Poincar\'e duality space $M$ of dimension $d$, the
induced map $q_X : \mbox{Map}(X, M) \to \mbox{Map}(S^1\amalg S^1,
M)$ is a fibration that is the pullback fibration of the usual
path fibration $M^{[0,1]} \to M\times M$ along the projection
$p\times p : LM\times LM \to M\times M$. We have recovered Diagram
$(\ast\ast)$ of the introduction:
$$
\xymatrix{
LM\times_MLM \simeq \mbox{Map}(X,M)\ar[d]_{\scriptstyle q_X} \ar[rr] && M^{[0,1]}
\ar[d]\ar[r]^{\simeq}& M \ar[ld]^{\Delta}& \ar@{}[d]^{.} \\
LM\times LM \ar[rr] ^{p\times p} && M\times
M &&
}$$
Denote by  $c: LM\times_MLM\to LM$ the composition of loops. The composite
$$ H^*(LM) \stackrel{H(c)}\to  H^*(LM\times_M LM)\stackrel{H({q_X}^!)}\to
 H^*(LM\times LM)$$ is, by definition,  the dual of the loop product $\bullet$ on
$H_*(LM)$.

\vspace{2mm} Since ${q_X}^!$ is a morphism of $C^*(LM\times
LM)$-modules,   $H_*({q_X}^!)  : H_*(LM\times LM)\to
H_*(LM\times_MLM)$ is   a morphism of $H^*(LM\times LM)$-modules
where the actions are given   by the cap product. Remark now that
the projection $p\times p : LM\times LM\to M\times M$ makes
$H_*(LM\times LM)$ a $H^*(M\times M)$-module.  We extend to the
case of Poincar\'e duality spaces   the following result  of
Tamanoi in (\cite{T2}, Theorem A),

\vspace{2mm}\noindent {\bf  Theorem 9.} {\sl Let $M$ be a
1-connected oriented Poincar\'e duality space, then the   loop
product is a morphism of $H^*(M)\otimes H^*(M)$-modules:
$$(\alpha_1 \cap b)\bullet (\alpha_2\cap c) = (-1)^{\vert
\alpha_2\vert (d+\vert b\vert)} (\alpha_1\cup\alpha_2)\cap
(b\bullet c)\,, \hspace{5mm} b,c\in H_*(LM),  \alpha_1, \alpha_2
\in H^*(M)\,.
$$}

\vspace{2mm}\noindent {\bf Proof.} First remark that $H_*({q_X}^!)
 : H_*(LM)\otimes H_*(LM) \to
H_*(LM\times_MLM)$ is a morphism of $H^*(LM)\otimes
H^*(LM)$-modules, and therefore  a morphism of $H^*(M)\otimes
H^*(M)$-module. This implies  that
$$H_*({q_X}^!) ((\alpha\cap b)\otimes ( \beta\cap c)) =
(-1)^{d(\vert \alpha\vert +\vert\beta\vert)+ \vert b\vert \cdot
\vert \beta\vert} (\alpha \cup \beta) \cap H_*({q_X}^!)((b\otimes
c)\,.$$ Now, since the composition of loops $LM\times_MLM\to LM$
is a morphism over $M$, the induced map in homology is a morphism
of $H^*(M)$-modules for the cap product. This implies the
result.\hfill $\square$

\subsection{The loop coproduct}

Let $M$ be a 1-connected oriented  Poincar\'e duality space. Consider the fibration $q :
Z \to LM$ where
$$Z=\{(\omega, c)\,\vert \,\omega\in LM, c\in M^{[0,1]}\,,
c(0)=\omega (0)\,, c(1)=\omega (1/2)\}$$ and $q(\omega,
c)=\omega$. The space $Z$ has the homotopy type of $LM\times_MLM$.
Denote by $\pi : Z\to LM\times LM$  the canonical injection. The
dual of the loop coproduct can be described as the composition of
$H^*(\pi)$ with the  shriek map $H^*(q^!)$,
$$H^*(LM)\otimes H^*(LM)\stackrel{H^*(\pi)}{\to}  H^*(Z)\stackrel{H^*(q^!)}{\to} H^*(LM)
\,.$$ The fibration $q$ is obtained as a pullback  of the
fibration $(p_0,p_1) : M^{[0,1]}\to M\times M$ along the map $\ell
: LM\to M\times M$ that evaluates a loop at the base point and at
the middle point, $\ell (\omega) = (\omega(0), \omega (1/2))$,
$$\xymatrix{
Z \ar@{^(->}[rr] \ar[d]_{\scriptstyle q}&& M^{[0,1]}\ar[d]^{\scriptstyle
(p_0,p_1)}&\ar@{}[d]^{.}\\
LM \ar[rr]^{\ell}&& M\times M&
}$$ Since $\ell$
is clearly homotopic to the composite $LM\stackrel{p}{\to} M
\stackrel{\Delta}{\to} M\times M$,  a model for $q^!$ is obtained
by pulling back a model for $(p_0,p_1)^!=\Delta^!$ along
$\Delta\circ p$.

\subsection{The intersection map with the fiber in a monoidal
fibration}

A monoidal fibration is a fibration $F \stackrel{j}{\to}
E\stackrel{p}\to B$, with a multiplication $\mu : E\times_BE\to E$
that extends a multiplication $\mu_0$ on $F$. The basic example is
the free loop space fibration. Denote by $\Delta': E\times_BE \to
E\times E$ the natural injection. Then the composition
$$H_*(E\times E)\stackrel{H_*(\Delta'^!) }{\to}
H_*(E\times_BE)\stackrel{H_*(\mu)}{\to} H_*(E)$$ defines a {\it $\mu$-intersection product} on $H_*(E)$. Generalizing the result of Chas and
Sullivan for the usual loop product, we have,

\vspace{2mm}\noindent {\bf Theorem 10.}  {\sl Suppose that $F,E,B$
are oriented Poincar\'e duality spaces and that $B$ is simply
connected. Then the intersection  map with the fiber $$H_*(j^!):
H_*(E) \to H_*(F)$$ is a multiplicative map  with respect to the
$\mu$-intersection product and the product $H_*(\mu_0)$.}

\vspace{2mm}\noindent {\bf Proof.} Denote by $k : \{b_0\} \to B$
the injection of the base point. Since the product $\mu_0$ is
the pullback of the multiplication $\mu$ along $k$, we have the
commutative diagram
$$
\xymatrix{
  &&&
F\times F \ar[rr]^{j'} \ar'[d][dd] \ar[dl]^{\mu_0} && E\times_BE
\ar[dd]  \ar[dl]^{\mu}&
\\
 && F \ar[rr]^(0.6){j } \ar[dr]  && E
\ar[dr] && \\
 &&&\{b_0\} \ar[rr]
&& B &
 }
$$
We deduce then from Theorem 4 that $H_*(\mu_0)\circ H_*(j'^!)  =
H_*(j^!) \circ H_*(\mu)\,.$ On the other hand  we consider  the
commutative diagram obtained by  obvious pullbacks
$$\xymatrix{
&F\times F \ar[rr]^{j\times j} \ar'[d][dd] && E\times E\ar[dd]^{p\times p}&\ar@{}[dd]^{.}
\\
F\times F \ar@{=}[ur] \ar[rr]^(0.3){j'} \ar[dd]&& E\times_BE
\ar[ur]^{\Delta'} \ar[dd]
\\
&\{ (b_0,b_0)\} \ar'[r][rr]&& B\times B&\\
\{ b_0\}\ar@{^(->}[rr]\ar[ur]^{\Delta}&&B \ar[ur]^{\Delta}\,. }
$$
By  uniqueness of the shriek map (Theorem 2) we obtain
$$H_*(j'^!) \circ H_*(\Delta'^!) = H_*(id^!) \circ H_*((j\times j)^!)  =
H_*(j^!) \otimes H_*(j^!) \,.$$ Then,
$$H_*(j^!)  \circ (H_*(\mu) \circ H_*(\Delta'^!) ) = H_*(\mu_0) \circ
H_*(j'^!) \circ H_*(\Delta'^!)  = H_*(\mu_0)\circ (H_*(j^!) \otimes
H_*(j^!) )\,.$$ \hfill $\square$

\section{Rational string topology}

\subsection{Rational homotopy theory}

To make computations over a field $\bk$ of characteristic zero,
the good tool is the theory of minimal models introduced by
Sullivan in \cite{Su} (see also \cite{FHT3}). For recall, Sullivan
defines a functor $A(-)$ from the category of topological spaces
to the category of commutative differential graded $\mathbb
Q$-algebras (for short cdga). One major property of the functor
$A(-)$ is the existence of a functor $E(-)$ with values in
differential graded algebras, and the existence of natural
quasi-isomorphisms
$$A(X) \leftarrow E(X) \to C^*(X;\mathbb Q)\,.$$

If $V$ is a graded vector space, $\land V$   denotes the free
commutative graded algebra on $V$. A {\sl Sullivan cdga} is a cdga
of the form $(\land V,d)$ such that $V$ admits a basis $(v_i)$
indexed by a well ordered set such that $d(v_i) \in \land (v_j,
j<i)$. A {\sl minimal cdga} is a Sullivan cdga $(\land V,d)$ for
which $d(V) \subset \land^{\geq 2}(V)$. Now if $(A,d)$ is a cdga
such that $H^0(A,d)=\mathbb Q$, there exists a minimal cdga
$(\land V,d)$ equipped with a quasi-isomorphism $\varphi : (\land
V,d) \to (A,d)$. This property characterizes $(\land V,d)$ up to
an isomorphism. The cdga $(\land V,d)$ is   called {\sl the
minimal model} of $(A,d)$. In particular the minimal model of
$A(X)$ is called the minimal model of $X$ and is denoted by ${\cal
M}_X$. More generally we call model of $X$ any  cdga which is
quasi-isomorphic to ${\cal M}_X$.

 Every
continuous map between connected spaces $f : X \to Y$ induces a
unique map  up to homotopy  ${\cal M}_f : {\cal M}_Y\to {\cal
M}_X$ making commutative up to homotopy the diagram
$$\xymatrix{
A(Y) \ar[rr]^{A(f)}&& A(X)&\ar@{}[d]^{.}\\
{\cal M}_Y\ar[u]^{\scriptstyle \simeq}  \ar[rr]^{{\cal M}_f}&&
{\cal M}_X \ar[u]_{\scriptstyle \simeq}& }
$$
 The  map ${\cal M}_f$ is called the {\it minimal
model of $f$}. More generally we call model of $f$ every map of
cdga's $h : (A,d) \to (B,d)$ such that there are
quasi-isomorphisms $\varphi$ and $\psi$ making commutative, up to
homotopy,  the diagram
$$
\xymatrix{ {\cal M}_Y \ar[rr]^{{\cal M}_f}\ar[d]_{\scriptstyle
\varphi}&& {\cal
M}_X \ar[d]^{\scriptstyle \psi}&\ar@{}[d]^{.}\\
(A,d) \ar[rr]^{h}&& (B,d)& }$$ Of course
  $\varphi$ and $\psi$ are parts of the structure
of the model $h$.

A {\it relative Sullivan model} for a homomorphism of cdga's  $f :
(A,d) \to (B,d)$ is an injection $(A,d) \hookrightarrow (A\otimes
\land V,d)$ equipped with a quasi-isomorphism $\varphi : (A\otimes
\land V,d) \to (B,d)$ satisfying $\varphi (a\otimes 1)=f(a)$, and
where $V$ admits a basis indexed by a well ordered set such that
$D(v_i) \subset A \otimes \land (v_j, j<i)$. The   relative model
is called {\it minimal} if $d(V) \subset \land^{\geq 2}(V)\oplus
(A^{\geq 1}\otimes \land V)$.

Let now $F\stackrel{j}{\to} E\stackrel{p}{\to} B$ be a fibration,
$(A,d)$ a model for $B$ and $\varphi : (A,d) \to (A\otimes \land
V,d)$ be a relative model for $p$. The morphism $\varphi$ is
called a relative model for the fibration $p$. The cdga $(\land
V,\bar d):=\mathbb Q\otimes_{A}(A\otimes \land V,d)$ is then a
model for   $F$.

 For instance a model  for the diagonal map $\Delta : X
\to X\times X$ is given by  the multiplication     of ${\cal M}_X
= (\land V,d) $, $\mu :  {\cal M}_X\otimes{\cal M}_X \to{\cal M}_X
$. A relative minimal model for the diagonal has the form
$$\rho : (\land V \otimes \land V\otimes \land sV,d)\to (\land V,d)$$
where $(sV)^n = V^{n+1}$ and $\rho (sV)=0$. The cdga $$(\land V
\otimes \land sV, D):= (\land V,d)\otimes_{\land V\otimes \land V}
(\land V\otimes \land V\otimes \land sV,d)$$ is then a minimal
model for the free loop space on $X$.  Then, $$D(sv) =-Sd(v)\,,$$
where $S$ is the   derivation defined by $S(v) = sv$ and $S(sv)=0$
(\cite{VS}, \cite[page 206]{FHT}).

  A
{\it Poincar\'e duality cdga of dimension $m$}  is a cdga $(A,d)$
satisfying:
\begin{enumerate}
\item[(i)] $A$ is finite dimensional, $A^{>m}=0$, $A^m = \mathbb Q \omega$,
 \item[(ii)] the map $\theta : A^r  \to
\mbox{Hom}(A^{m-r}, \mathbb Q)$  defined by $ab = \theta(a)(b)
\cdot \omega $ is an isomorphism for $0\leq r\leq m$.
\end{enumerate}
 By a result
of Lambrechts and Stanley \cite{LS} any  simply connected space
whose rational cohomology satisfies Poincar\'e duality admits a
Poincar\'e duality model.

\subsection{The intersection map with the fiber with rational coefficients}

Here we will use Sullivan models to give  a  rational description
of $H^*(j^!)$ for a fibration $ F \stackrel{j} \to E \to M$  with
base  a 1-connected oriented Poincar\'e duality space.

  Let $(A,d)\to (A\otimes \land V,D) \to (\land V,
\overline{D})$ be a relative Sullivan model for the fibration $p :
E\to M$. We suppose that  $A^0=\mathbb Q$ and $A^{>m}=0$.  If
$\omega\in A^m$ is a cocycle representing the fundamental class of
$M$, then a model of $i^!$ is given by
$$i^! : \bk \to A \,, \hspace{1cm} i^!(1) = \omega\,.$$
Therefore, following the proof of Theorem 7 we obtain:

\vspace{2mm}\noindent {\bf Proposition 1.} {\sl The map $H^*(j^!)
: H^*(F;\mathbb Q)\to H^*(E;\mathbb Q)$ is induced by the morphism
of $(A\otimes \land V,D)$-modules
$$  (\land V,\overline{D}) \to (A\otimes \land V,D)\,, \hspace{1cm}
 \alpha \mapsto (-1)^{\vert \omega\vert \cdot \vert
\alpha\vert} \omega \otimes \alpha\,.$$}

 \vspace{2mm}\subsection{The rational loop product}

Let $M$ be a Poincar\'e duality space of dimension $m$ and let
$(A,d)$ be a Poincar\'e duality model for $M$. Denote by $(a_i)_i$
a homogeneous  basis of $A$ and by $(a'_i)_i$ its Poincar\'e  dual
basis:  $a_i\cdot a'_j = \delta_{ij}\,\omega$. The diagonal
element
$$ {D} = \sum_i (-1)^{\vert a_i\vert} a_i\otimes a'_i  \in A \otimes A$$
  is a cocycle such that,  for each element $a\in A$,  $(a\otimes 1) {D} =
  (1\otimes a) {D}$.
  The multiplication by
$ {D}$ is thus a non-trivial morphism of $(A\otimes A)$-modules
$$\mu_{ {D}} \colon A \to A\otimes A\,.$$
Since $\mu_{ {D}}$ maps the fundamental class to the fundamental
class, by Theorem 1, $\mu_D$ is a representative for $\Delta^!$
where   $\Delta : M \to M^2$ denotes  the diagonal map.

\vspace{2mm}
   Denote by $(A\otimes \land sV,D)$ be a Sullivan   model
   for the free loop space $LM$.
As a corollary of Theorem 2,  we obtain

\vspace{2mm}\noindent {\bf Proposition  2.}  {\sl  With the above
notation, the morphism
$$\mu_{ {D}}\otimes 1 : A\otimes_{A\otimes A}(A\otimes \land V,D)^{\otimes
2}\to (A\otimes A)\otimes_{A\otimes A} (A\otimes \land V)^{\otimes
2}= (A\otimes \land V)^{\otimes 2}$$ induces in cohomology the
  map $H^*(q^!): H^*(LM\times_MLM) \to
H^*(LM\times LM)$.}

\vspace{2mm} This allow us to extend to Poincar\'e duality space
the next result of \cite{FT}:

\vspace{2mm}\noindent {\bf Theorem 11.}     {\sl There exists  an
isomorphism of  graded algebras between the Hochschild cohomology
of $A$, $H\!H^*(A)$, and the loop homology $H_*(LM)$. This
isomorphism identifies in cohomology the dual of the loop product
with the map induced by the composite
$$C\!H_*(A) \stackrel{\nabla}{\to}
A\otimes_{A^{\otimes^2}} (CH_*(A))^{\otimes^2}
\stackrel{\mu_D\otimes 1}{\to} A^{\otimes^2}
\otimes_{A^{\otimes^2}} (CH_*(A))^{\otimes^2}=
(CH_*(A))^{\otimes^2}\,.$$ \noindent where  $\nabla : C\!H_*(A)\to
A\otimes_{A^{\otimes^2}} (C\!H_*(A))^{\otimes^2}$ is the morphism
of complexes defined by
$$\nabla \left(a\otimes [a_1\vert \cdots
\vert a_n]\right)= \sum_{i=0}^n \, a\otimes [a_1\vert \cdots \vert
a_i] \otimes [a_{i+1}\vert \cdots \vert a_n]\,.$$ }

\vspace{2mm}\noindent {\bf Proof.} As proved in \cite{FT2}, the
map $\nabla$ is a semifree model for the composition of loops $c :
LM\times_MLM \to LM$. On the other hand, denoting by $(A\otimes
\land V,D)$ the model for the free loop space described above,
  there is a quasi-isomorphism of semifree $A$-modules $\psi :
(A\otimes \land V,D) \to C\!H_*(A)$.

Denote by $m$ the multiplication in $A$ and by $\theta : A \to
A^\#=\mbox{Hom}(A,\mathbb Q)$ a representative of the cap product
with the orientation class. Then  we obtain the commutative
diagram whose vertical lines are quasi-isomorphisms
$$
\xymatrix{ A \ar[rr]^{\mu_D} \ar[d]_{\scriptstyle \theta}
&&A\otimes A\ar[d]^{\scriptstyle
 \theta\otimes\theta}& \\
A^\# \ar[rr]^{m^\vee}&& A^\# \otimes A^\#\,.}$$ This is the main
fact which allows to   end the proof as  in (\cite{FT2}). \hfill
$\square$

\subsection{The rational loop coproduct}

 Suppose $M$ is a 1-connected oriented  Poincar\'e
duality space and $\bk=\mathbb Q$. We use the notations of section
3.3.  To obtain a model for the loop coproduct we need first to
obtain a model for $q^!$ in the derived category of
$C^*(LM)$-modules. This model is
 obtained from a model of $(p_0,p_1)^!$ by tensorization  with a
model of the map $\ell$.   We will recover in this way the model
described  in \cite{CT}.

Let $(A,d)$ be a Poincar\'e duality model for $M$ with fundamental
class $\Omega$, and $\theta : (A\otimes A\otimes \land Z,d)\to
(A,d)$,
$$\theta (a\otimes
a'\otimes 1)=aa'\hspace{5mm}\mbox{and}\hspace{5mm} \theta (Z)=0 $$
be a relative Sullivan  model for the product $A\otimes A\to A$.
   Let $(A\otimes
A\otimes \land  {Z'}, D)$ be a copy of $(A\otimes A\otimes \land
 {Z}, D)$ and form the tensor product
$$(A\otimes A\otimes \land {Z}\otimes \land
 {Z'}, D):= (A\otimes A\otimes \land  {Z},
D)\otimes_{A\otimes A}(A\otimes A\otimes \land  {Z'}, D)\,.$$ Then
the injection $$ \rho : (A,d)\otimes (A,d) \to (A\otimes A\otimes
\land  {Z}\otimes \land  {Z'}, D)$$ is a model for $\ell : LM\to
M\times M$.   Recall now that the shriek map associated to
$(p_0,p_1)$ is represented by the multiplication by the diagonal
class $\mu_{ {D}} : A\to A\otimes A$ (see section 4.1). Therefore
a model for the shriek map associated to $q$ is given by $q^!$ in
the diagram
$$\xymatrix{
(A \otimes \land  {Z}\otimes\land  {Z'},D) \ar@{=}[rr]
\ar[d]_{\scriptstyle q^!}&& A\otimes_{A\otimes A} (A\otimes A
\otimes \land  {Z}\otimes\land  {Z'},D) \ar[d]^{\scriptstyle \mu_{
{D}}\otimes 1}\\
(A \otimes A \otimes \land  {Z}\otimes\land  {Z'},D)\ar@{=}[rr]&&
(A\otimes A)\otimes_{A\otimes A}(A\otimes A\otimes \land
 {Z}\otimes\land  {Z'},D)
 }$$
 Now observe that the projection $\theta \otimes 1 :
(A\otimes A\otimes \land \ {Z}\otimes\land  {Z'},D) \to (A \otimes
\land  {Z'},D)$ is a quasi-isomorphism, and form the composition
$$\psi = (\theta\otimes 1)\circ q^! : (A \otimes   \land  {Z}\otimes\land  {Z'},D) \to
(A\otimes \land  {Z'},D)\,.$$ We have $$\psi (a\otimes b\otimes c)
= \left\{ \begin{array}{l} 0\,, \hspace{3mm}\mbox{if $a\otimes
b\in
(A\otimes \land Z)^+$}\\
 \chi (M)\cdot\Omega   \otimes c\,, \hspace{3mm}
  \mbox{if $a\otimes b=1\otimes 1$}\end{array}\right.$$
  Here $\chi
(M)$ denotes the  Euler-Poincar\'e characteristic of $M$.
Therefore, if $\chi (M)=0$ then  $H^*(q^!)=H^*(\psi) = 0$,  and
thus the loop coproduct is trivial.

\section{Gorenstein spaces}

Let $X$ be  simply connected $\bk$-Gorenstein space of dimension
$d$ whose cohomology $H^*(X;\bk)$ is of finite type. The diagonal
map $\Delta : X \to X^n$ makes $C^*(X)$ into a $C^*(X^n$)-module,
and we have

\vspace{2mm}\noindent {\bf Theorem 12.}
$$\mbox{Ext}_{C^*(X^n)}(C^*(X), C^*(X^n)) \cong s^{(n-1)d}\,
H^*(X)\,.$$

 \vspace{2mm}\noindent {\bf
Proof.} The proof will proceed in two main steps, the first one
consisting to replace the statement concerning the cochains by a
statement concerning their free models (§1.4). The second one will concern
the computation of the corresponding Ext.

\vspace{2mm} \noindent {\bf Step 1 of the proof.} (Construction of
models). Let
$$
\varphi_1 : (T(V),d) \stackrel{\simeq}{\to} C^*(X;\bk) \mbox{ and }
\varphi_{n-1}: (T(W),d) \stackrel{\simeq}{\to} C^*(X^{n-1};\bk)
$$
be free models of $X$ and $X^{n-1}$.

The tensor product $(T(V),d)\otimes (T(W),d)$ admits a free model
of the form
$$\varphi :  (T(V\oplus W \oplus s(V\otimes W)),D) \to
(T(V),d)\otimes (T(W),d)
$$
  where the differential $D$ extends the   differentials $d$
already defined on $T(V)$ and $T(W)$, and  $\varphi(v)=v$,
$\varphi(w)=w$ and $\varphi (s(V\otimes W))=0$. For sake of
simplicity we write  $ V\lozenge W :=V\oplus W \oplus s(V\otimes
W)$.

Consider the commutative diagram

\centerline{ \footnotesize $
\xymatrix{
\ar@{}[d]^{(\ast\ast\ast)} &C^*(X) \ar[rr]^{C^*(q_1)}\ar[d]^{\scriptstyle 1\otimes id} && C^*(X^n)\ar[d]^{\scriptstyle
(EZ)^\#}&\ar@{}[d]^{.}\\
&C^*(X^{n-1})\otimes C^*(X) \ar[rr]^{\simeq}&&
(C_*(X^{n-1})\otimes C_*(X))^\#&
}$}

\vspace{1mm}\noindent in which  $EZ$  denotes the Eilenberg-Zilber
chain equivalence and  $q_1$ denotes the projection on the last
factor. Since $(EZ)^\#$ is surjective, by (FE2)  and the diagram
$(\ast\ast\ast)$ above the map
$$\varphi_n := C^*(q_1)\circ \varphi_1 : (T(V),d) \to
C^*(X^n)
$$
 extends into a morphism of differential graded algebras,
also denoted $\varphi_n$,  that makes commutative the next
diagram.

\centerline{\footnotesize $
\xymatrix{
(T(V\lozenge W),D) \ar@{.>}[rrrr]^{\varphi_n} \ar[d]_{\scriptstyle \simeq}
&&&&  C^*(X^n) \ar[dd]^{\scriptstyle EZ^{\#} } _{\simeq}&\ar@{}[dd]^{.} \\
(T(W),d)\otimes (T(V),d) \ar[d]_{\scriptstyle \varphi_{n-1}
\otimes\varphi_1} ^{\simeq}&& (T(V),d) \ar@{_(->}[llu] \ar[rru]^{\varphi_n}\\
 C^*(X^{n-1})\otimes C^*( X) \ar[rrrr]^{\simeq }
  &&&& \left(C_*(X^{n-1})\otimes C_*(X)\right)^{\#}&
}$}

\vspace{1mm}\noindent  Choose  $e\in X$, and  denote by $i_{n-1} :
X^{n-1} \to X^n$ and $j_1 : X\to X^n$ the injections defined by
$i_{n-1}(a)=(a,e)$, and $j_1(a) = (e,\ldots ,e,a)$. Denote also by
$\varepsilon : C^*(X) \to C^*(\{ e\}) $ the augmentation map.
 Since we work with normalized cochain complexes,   $C^*(\{ e\})=
\bk$.

  The commutativity of the diagram

\centerline{\footnotesize $\xymatrix{
C^*(X^n)  \ar[d]_{\scriptstyle EZ^\#}^{\simeq} \ar[drrr]^{C^*(i_{n-1})} \\
(C_*(X^{n-1})\otimes C_*(X))^\# \ar[rrr]^{C_*(i_{n-1})^\#}  &&& C^*(X^{n-1})\\
C^*(X^{n-1})\otimes C^*(X) \ar[u]  \ar[urrr]^{1\otimes
\varepsilon} }$}

\vspace{1mm}\noindent implies  the commutativity of the diagram of
differential graded algebra,

$$\xymatrix{
(T(V\lozenge W),D) \ar[rr]^{\varphi_n}_{\simeq} \ar[d]_{\scriptstyle (1\otimes\varepsilon)\circ \varphi}&& C^*(X^n)\ar[d]^{\scriptstyle
C^*(i_{n-1})}\\
(T(W),d) \ar[rr]^{\varphi_{n-1}}_{\simeq}&& C^*(X^{n-1})}$$
 Therefore
$$C^*(i_{n-1})\circ \varphi_n ( V\oplus s(V\otimes W))=0\,.
$$

We   now apply to $C^*(\Delta)$ the lifting homotopy property for
free models.  Since $q_1\circ \Delta = id_X$ and
$\varphi_n(v)=C^*(q_1)\varphi_1(v)$ the identity on $T(V)$ extends
into a morphism $\psi  $ making commutative
  the following diagram of differential graded algebras
$$
\xymatrix{
 (T(V\lozenge W), D) \ar[rr]_{\simeq}^{\varphi_n } \ar[d]_{\scriptstyle \psi}&&
C^*(X^n)\ar[d]^{\scriptstyle
C^*(\Delta)}\\
  (T(V),d) \ar[rr]^{\varphi_1}_{\simeq}&& C^*(X)
}
$$
 The morphism $\psi$ makes $T(V)$ a $T(V\lozenge W)$-module. By construction
the induced structure of $T(V)$-module on the submodule  $T(V)$ of
$T(V\lozenge W)$ coincide with the  usual multiplication on
$T(V)$.  We  have the isomorphisms
$$\renewcommand{\arraystretch}{1.6}\begin{array}{l}\mbox{Ext}_{C^*(X^n)}(C^*(X), C^*(X^n))
\cong
\mbox{Ext}_{T(V\lozenge W)}(T(V), T(V\lozenge W))\\
\cong \mbox{Ext}_{T(V\lozenge W)}(T(V), T(W)\otimes T(V))\\
\cong H^*\left(\,\mbox{Hom}_{T(V\lozenge W)}( \,(T(V \lozenge
W)\otimes Z,D), T(W)\otimes
T(V)\,)\,\right)\,.\end{array}\renewcommand{\arraystretch}{1}$$
\noindent Here $(T(V\lozenge W)\otimes Z,D)$ denotes a
$T(V\lozenge W)$-semifree model of $T(V)$.

\vspace{3mm}\noindent {\bf Step 2 of the proof.}
(Computation of $\mbox{Ext}_{C^*(X^n)}(C^*(X), C^*(X^n))$). \\
We fix a non negative integer $N$ and we  filter the complex
$$\mbox{Hom}_{T(V\lozenge W)}\left( T(V\lozenge W)\otimes Z,\,
T(W)\otimes \frac{T(V)}{T(V)^{>N}}\right)$$ by the sub vector
spaces, $F^p \subset F^{p-1}$ defined  by:
$$F^p =  \{ \, \varphi \, \vert \, \varphi (Z) \subset T(W)
\otimes \frac{T(V)^{\geq p}}{T(V)^{>N}}\, \}\,.$$
  This filtration  induces a
converging spectral sequence with
$$E_0^{p,*} = \{\, \varphi\,\vert\, \varphi(Z)\subset T(W) \otimes
T(V)^p\,\}\,,$$

\vspace{2mm}\noindent {\bf Lemma 2.} {\sl This  spectral sequence
collapses at the $E_2$-level and
$$E_2^{p,*} \cong s^{(n-1)d}
H^p(T(V)/T(V)^{>N})\,.$$}

\vspace{2mm}\noindent {\bf Proof.}   Recall that
$V\lozenge W = V \oplus W \oplus s(V\otimes W)$ and denote by $I$ the ideal in
$T(V \lozenge W)$ generated by $V$ and $s(V\otimes W)$.   Since $I$ acts
trivially on $T(W) \otimes T(V)^p$, we have an isomorphism of
complexes
$$\renewcommand{\arraystretch}{1.4}
\begin{array}{l}\mbox{Hom}_{T(V\lozenge W)} (T(V\lozenge W)\otimes Z,
T(W)\otimes T(V)^p) \\
\simeq \mbox{Hom}_{T(W)}(T(W)\otimes_{T(V\lozenge W)}(T(V\lozenge
W)\otimes Z),\, T(W) \otimes T (V)^p)\,.\end{array}
\renewcommand{\arraystretch}{1}$$
 Here $T(W)$ is
considered as a $T(V\lozenge W)$-module via the morphism
$(\varepsilon\otimes 1)\circ \varphi$ . We have also
quasi-isomorphisms between bar constructions:
$$
\renewcommand{\arraystretch}{1.3}
\begin{array}{ll}
\mathbb B(T(W), T(V\lozenge W), \bk) &\simeq \mathbb B(T(W),
T(W)\otimes T(V), \bk)\\ &\cong    \mathbb B (T(W), T(W),
\bk)\otimes \overline{\mathbb B}(T(V)) \simeq \overline{\mathbb
B}(T(V))\,.\end{array}\renewcommand{\arraystretch}{1}
$$
 It follows that the natural injection $$\bar \mathbb B( T(V)):=\mathbb B(\bk, T(V),\bk)
\hookrightarrow \mathbb B(T(W), T(V\lozenge W), \bk)$$ is a
quasi-isomorphism of complexes. We consider then the injection
$$
j : \mathbb B(\bk, T(V), T(V\lozenge W)) \hookrightarrow \mathbb
B(T(W), T(V\lozenge W), T(V\lozenge W))\,,$$ and we filter the
complexes $\mathbb B(\bk, T(V), T(V\lozenge W))$ and $\mathbb
B(T(W), T(V\lozenge W), T(V \lozenge W))$ respectively by
$\overline{\mathbb B}(T(V))^{\geq p} \otimes T(V\lozenge W)$ and
$\left( T(W)\otimes \overline{\mathbb B}(T(V\lozenge
W))\right)^{\geq p} \otimes T(V \lozenge W)$. The morphism $j$
preserves the filtrations and induces an isomorphism at the
$E_2$-level of the associated spectral sequences. Therefore $j$ is
a quasi-isomorphism so that we deduce the isomorphisms
$$\renewcommand{\arraystretch}{1.6}\begin{array}{l}
H_*(T(W)\otimes_{T(V\lozenge W)}(T(V\lozenge W)\otimes Z)) \cong
\mbox{Tor}_*^{T(V\lozenge W)}(T(W),
T(V))\\
\cong H_*( \mathbb B (T(W), T(V\lozenge W), T(V\lozenge W)) \otimes_{T(V\lozenge W)} T(V))
\\
\cong H_*(\mathbb B(\bk, T(V), T(V\lozenge W))
\otimes_{T(V\lozenge W)}T(V)) \cong H_*(\mathbb B(\bk, T(V),
T(V))) \cong \bk\,.\end{array}\renewcommand{\arraystretch}{1}$$
 This shows that $T(W)\otimes_{T(V\lozenge
W)}(T(V\lozenge W)\otimes Z)$ is a semifree resolution of $\bk$ as
$T(W)$-module. Therefore,
$$E_1^{p,*}  \cong
\mbox{Ext}_{T(W)}(\bk,T(W))\otimes T(V)^{p}\cong
s^{(n-1)d}T(V)^p/T(V)^{>N} $$ and $$E_2  \cong s^{(n-1)d}
H^*(T(V)/T(V)^{>N}) \,.$$

\hfill $\square$

\vspace{3mm}\noindent  Now since $T(V)$ is of finite type,
$T(W)\otimes T(V)=\displaystyle \lim_{\leftarrow_N}\, T(W)\otimes
\frac{T(V)}{T(V)^{>N}}$, and
$$\renewcommand{\arraystretch}{1.2}\begin{array}{l}
\mbox{Hom}_{T(V\lozenge W)}\left(T(V\lozenge W)\otimes
Z,T(W)\otimes T(V)\right) \\ = \displaystyle\lim_{\leftarrow_N}
\mbox{Hom}_{T(V\lozenge W)} \left( (T(V \lozenge W)\otimes Z, T(W)
\otimes\frac{T(V)}{T(V)^{>N}}\right)
\,.\end{array}\renewcommand{\arraystretch}{1}$$ We get therefore
the short exact sequence
$$\renewcommand{\arraystretch}{1.6}
\begin{array}{ll} 0 &\to {\displaystyle\lim_{\leftarrow_N}}^1 s^{(n-1)d}
H^{p+1-(n-1)d}(T(V)/T(V)^{>N})  \to \mbox{Ext}^p_{T(V\lozenge W)}
(T(V),
T(W)\otimes T(V)) \\
& \to
\displaystyle\lim_{\leftarrow_N}s^{(n-1)d}H^{p-(n-1)d}(T(V)/T(V)^{>N})\to
0\end{array} \renewcommand{\arraystretch}{1}
$$
 Since the tower
 $ (s^{(n-1)d}H^p(T(V)/T(V)^{>N}) )_N $  satisfies the  Mittag-Leffler condition,
 we  have in fact    isomorphisms
$$\mbox{Ext}^p_{T(V\lozenge W)}(T(V), T(W) \otimes T(V)) \cong \lim_{\leftarrow_N}
s^{(n-1)d} H^p(T(V)/T(V)^{>N})\cong H^{p-(n-1)d}(X)\,.$$

\hfill $\square$

 \vspace{5mm}\noindent {\bf Corollary.} {\sl Let
$\Delta : X^r \to X^n$ be the product of   diagonal maps $
X \to X^{n_i}$, $i=1, \ldots , r$. Then
$$\mbox{Ext}_{C^*(X^n)}(C^*(X^r), C^*(X^n)) \cong s^{(n-r)d}\,
H^*(X^r)\,.$$  Here $C^*(X^r)$ is viewed as a $C^*(X^n)$-module
via $C^*(\Delta)$. }

\vspace{2mm} \noindent {\bf Proof.}  Write $\Delta$ as the
product of diagonal maps $$\Delta = \Delta_1 \times \cdots \times
\Delta_n : X \times \cdots \times X \to X^{n_1}\times \cdots
\times X^{n_r}\,.$$ Then denoting by $T_{X^s}$ a free model for
$C^*(X^s)$, we have

 \centerline{\footnotesize$\renewcommand{\arraystretch}{1.6}\begin{array}{l} \mbox{Ext}_{C^*(X^n)} (C^*(X^r), C^*(X^n)) \cong
\mbox{Ext}_{T_{X^n}}(T_{X^r}, T_{X^n})\\
\cong \mbox{Ext}_{T_{X^n}}(\otimes_{i=1}^r T_X, T_{X^n})\cong
\mbox{Ext}_{\otimes_{i=1}^r T_{X^{n_i}}} (T_X, \otimes_{i=1}^r
T_{X^{n_i}})\cong  \otimes_{i=1}^r \mbox{Ext}_{T_{X^{n_i}}}(T_X,
T_{X^{n_i}})\\
\cong \otimes_{i=1}^r \mbox{Ext}_{C^*(X^{n_i})}(C^*(X),
C^*(X^{n_i})) \cong s^{(n-r)d} H^*(X^r)
\end{array}\renewcommand{\arraystretch}{1}$}
\hfill $\square$

\vspace{2mm} We denote by $\Delta^!$ the map defined in the
derived category of $C^*(X^n)$-module from $C^*(X^r) $ to
$C^*(X^n)$ corresponding to a generator of
$\mbox{Ext}^{(n-r)d}_{C^*(X^n)}(C^*(X^r), C^*(X^n))$. This element
is well defined up to homotopy and up to the multiplication by a
scalar.

Let  $p : E \to X^n$ be a fibration and consider the homotopy
pullback

\centerline{\footnotesize$\xymatrix{
E' \ar[rr]^{g} \ar[d]_{\scriptstyle p'}&& E \ar[d]^{\scriptstyle p}\\
X^r \ar[rr]^{\Delta} && X^n
}
$}

\vspace{1mm}\noindent In the same way that Theorem 2  and  Theorem
3 are proved, we get

\vspace{2mm}\noindent {\bf Theorem 13.} {\sl If the spaces arising in the above
diagram are Gorenstein spaces then  in the
derived category of $C^*(E)$-modules, there exists a unique morphism (up to
homotopy and up to multiplication by a scalar)  $g^! : C^*(E') \to C^*(E)$ making
commutative the
diagram

\centerline{\footnotesize$\xymatrix{
C^*(E') \ar[rr]^{g^!}&& C^*(E)& \\
C^*(X^r) \ar[u]^{\scriptstyle C^*(p')} \ar[rr]^{\Delta^!} &&
C^*(X^n)\ar[u] _{\scriptstyle C^*(p)}\,.& }$}

\vspace{1mm}\noindent When  the fiber of $\Delta$ is a Poincar\'e
duality space, then $H^*(\Delta^!)$ and $H^*(g^!)$ coincide, up to
multiplication by a scalar, with the integration along the fiber.
}

\section{String operations on a Gorenstein space}

\vspace{2mm} Let $S$ be a connected surface of genus $g$ with $p\geq 1 $ incoming
boundary components and $q\geq 1$ outgoing boundary components. Then the
injection of the incoming components defines for each Gorenstein
space $X$ a map $q_S : \mbox{Map}(S,X)\to (LX)^p$ that is the
homotopy pullbak in a diagram
$$\xymatrix{ \mbox{Map}(S,X) \ar[rr]^{q_S} \ar[d] &&(LX)^p
\ar[d]^{\scriptstyle \psi}\\
X^r  \ar[rr]^{\Delta} && X^t\,.}
$$
 By Theorem 12,
 in the derived category of $C^*((LX)^p)$-modules there exists  a well
defined map  (up to  homotopy and up to  multiplication by a scalar)
$$(q_S)^! : C^*(\mbox{Map}(S,X) )
\to C^*((LX)^p)\,.
$$ Thus, as descibed in the introduction in case of Poincaré duality spaces,  $(g, p+q)$-string operations are well defined up to  multiplication by a scalar.

\subsection{The loop product for a classifying space}

\vspace{2mm} First remark that if $X  $ is  the
classifying space $BG$ of a compact connected Lie group $G$ then
the homotopy fiber of the diagonal map $\Delta : X^r \to X^t$ is
the space $(\Omega BG)^{t-r} \cong G^{t-r}$ that is a finite
dimensional Poincar\'e space. Therefore, by Theorem 3, the
integration along the fiber gives a shriek map $\Delta^!$ as required by  Chataur
and Menichi,  \cite{CM}.

 \vspace{2mm} \noindent {\bf Theorem 14.}  {\sl If $\bk=\mathbb Q$ and  $G$ is a compact
 connected Lie group   then  loop product $H_*(LBG)\otimes H^\ast(LBG) \to H_*(LBG)
$ is  trivial.}

\vspace{2mm}\noindent {\bf Proof.}  A Sullivan  model for
$BG$ is given by $A_G = (\land (x_1, \ldots , x_n),0)$ where $n=$rank$G$ and the
$x_i$ have even degree,  $\vert x_i\vert = 2n_i$. A $A_G\otimes A_G$-semifre model of
$A_G$  is:
$$\varphi : (\land (x_1, \ldots , x_n, x_1', \ldots, x_n', \overline{x_1},
\ldots , \overline{x_n}), D) \to (\land (x_1, \ldots ,
x_n),0)\,,
$$
\noindent  where
$\vert x_i'\vert = 2n_i$ , $\vert
\overline{x_i}\vert =2n_i-1$, $D(x_i) = D(x_i')=0$,  the differential on the lefthand cdga is determined  by:
$D(\overline{x_i})= x_i-x_i'$ and the morphism of cdga's $\varphi$ is defined by :  $\varphi (x_i) = \varphi(x_i') =
x_i$ and $\varphi (\overline{x_i}) = 0$. Then, a representative of $\Delta^!$ is given by the morphism of $(\land
(x_i, x_i'),0)$-modules
$$\Delta^! : (\land (x_i, x_i', \overline{x_i}),D) \to (\land (x_i,
x_i'),0)$$ \noindent
 defined by $\Delta^! (\overline{x_{i_1}}\,
\overline{x_{i_2}}\cdots \overline{x_{i_r}}) = 0$ if
$\overline{x_{i_1}}\, \overline{x_{i_2}}\cdots
\overline{x_{i_r}}\neq \overline{x_1}\cdots \overline{x_n}$ and is
equal to $1$ otherwise. Since the differential in the minimal model of $BG$ is zero, the
Sullivan minimal  model for  the free loop space on $BG$ is
$(\land (x_1, \ldots , x_n, \widehat{x_1}, \ldots ,
\widehat{x_n}), 0)$ with $\vert \widehat{x_i} \vert= 2n_i-1$. A
model for the projection $p\times p : LM\times LM\to M\times M$ is
given by the inclusion $\rho : (\land (x_i, x_i'),0) \to (\land
(x_i, x_i', \widehat{x_i}, \widehat{x_i'}),0)$. A representative
for $(q_X)^!$ is thus given by the tensor product
$$\psi = (\land (x_i, x_i', \widehat{x_i}, \widehat{x_i'}),0)\otimes_{\land (x_i, x_i')}
\Delta^!   : (\land (x_i, x_i', \overline{x_i}, \widehat{x_i},
\widehat{x_i'}), D)\to (\land (x_i, x_i', \widehat{x_i},
\widehat{x_i'}), 0)\,.$$ \noindent  Now remark that  the injection
of $(\land (x_i, \widehat{x_i}, \widehat{x_i'}),0)$ into the first
factor is a quasi-isomorphism and the composition with $\psi$ is
zero. This shows that the loop product is zero for $BG$. \hfill
$\square$

\vspace{2mm} \noindent {\bf Remark.} Denote by $M_1\subset M_2 \cdots \subset M_n \subset
BG$ be a sequence of compact manifolds in $BG$ such that $BG =
\cup M_n$. Then by restriction to $M_n$ we get a fibration $\Omega
BG \to L_nBG \to M_n$ and a loop product on $H_*L_nBG$. The
injection $M_{n-1}\to M_n$ induces a multiplicative shriek map
$H_*(L_{n}BG) \to H_*(L_{n-1}BG)$. In \cite{GS} Gruher and
Salvatore define the loop product on $BG$ to be the  graded
algebra $\displaystyle\lim_{\leftarrow} H_*(L_nBG)$. This  so-called
loop product  is no more a product on $H_*(BG)$ and is different
from our loop product.

\subsection{The loop coproduct of a classifying space}

\vspace{2mm}\noindent {\bf Theorem 15.}  {\sl If $\bk=\mathbb Q$ and if $G$  is a compact Lie group  then  the loop coproduct $H_*(LBG;\mathbb Q) \rightarrow
H_*(LBG)\otimes H^\ast(LBG)$ is an injective map.}

\vspace{2mm}\noindent {\bf Proof.}  Denote by $(\land (x_i), 0)$ a Sullivan
minimal model for $BG$ as in the previous section. We consider anew
the first diagram   in § 3.4, with $X$ in place of $M$. We use the model of $\Delta^!$ described above,
$$\Delta^! : (\land (x_i, x_i', \overline{x_i}, D) \to (\land
(x_i, x_i'), 0)\,,$$   where $D(\overline{x_i})= x_i-x_i'$. A
model for $\ell$ is given by the inclusion $$\varphi : (\land
(x_i, x_i'), 0) \to (\land (x_i, x_i', \widetilde{x_i},
\widehat{x_i}), D)$$ where $\vert \widetilde{x_i}\vert = \vert
\widehat{x_i}\vert = 2n_i-1$, $D(\widetilde{x_i})= x_i-x_i'$ and
$D(\widehat{x_i}) = 0$. Therefore a model for $q^!$ can be choosen
as
$$q^!= \Delta^!\otimes 1 :   (\land (x_i, x_i', \overline{x_i}, \widehat{x_i},
\widetilde{x_i}),D) \to (\land (x_i, x_i',  \widehat{x_i},
\widetilde{x_i}),D)\,.
$$  We consider now the following diagram
where the vertical maps are quasi-isomorphisms of cdga's
$$\xymatrix{
(\land (x_i, x_i', \overline{x_i}, \widehat{x_i},
\widetilde{x_i}),D) \ar[rr]^{q^!}&& (\land (x_i, x_i',
\widehat{x_i},
\widetilde{x_i}),D)\ar[d]_ {\scriptstyle \psi}^{\simeq}& \\
(\land (x_i,   \overline{x_i}, \widehat{x_i}, ),0)
\ar[u]^{\scriptstyle \tau}_{\simeq} && (\land
(x_i,
  \widehat{x_i} ),D)\,.&
}$$ where the quasi-isomorphisms $\tau$ and $\psi$ are defined by:
$\tau (x_i) = x_i$, $\tau (\overline{x_i})=
\overline{x_i}-\widetilde{x_i}$, $\tau (\widehat{x_i})=
\widehat{x_i}$, $\psi (x_i) = \psi (x_i') = x_i$, $\psi
(\widetilde{x_i})= 0$ and $\psi (\widehat{x_i})= \widehat{x_i}$.
The composition $\psi \circ q^!\circ \tau$ is a surjective map
and so is $H(q^!)$.

A model for the injection
$LBG\times_{BG}LBG\to LBG\times LBG$ is given by the surjection
$$(\land (x_i, \overline{x_i}),0) \otimes (\land (x_i',
\widehat{x_i}),0)\to (\land (x_i, \overline{x_i},
\overline{x_i}),0)\,.$$ The dual of the loop coproduct is the
composition
$$H^*(LBG \times LBG)\stackrel{H^*(\pi)}{\to} H^*(Z) \stackrel{H^*(q^!)}{\to}
H^*(LBG)\,.$$  Since both maps are surjective, this is also true
for the composition.\hfill $\square$

 \vspace{1cm}  Universit\'e Catholique de Louvain,
1348, Louvain-La-Neuve, Belgium

 Universit\'e d'Angers, UMR CNRS- 49045 Bd Lavoisier, Angers,
France

\end{document}